\documentclass[11pt]{article}
 \usepackage{amsfonts}
 \usepackage{mathrsfs}
 \usepackage{bbm}
 \usepackage{amsthm}
\usepackage{amsmath,amssymb}
\def\MatrixFont{\bf}
\def\VectorFont{\bf}

\newcommand{\mB}{{\MatrixFont B}}

\newcommand{\mE}{{\MatrixFont E}}
\newcommand{\mF}{{\MatrixFont F}}

\newcommand{\mH}{{\MatrixFont H}}
\newcommand{\mI}{{\MatrixFont I}}

\newcommand{\mP}{{\MatrixFont P}}
\newcommand{\mQ}{{\MatrixFont Q}}
\newcommand{\mR}{{\MatrixFont R}}

\newcommand{\mX}{{\MatrixFont X}}

\newcommand{\vb}{{\VectorFont b}}

\newcommand{\vu}{{\VectorFont u}}
\newcommand{\vv}{{\VectorFont v}}
\newcommand{\vw}{{\VectorFont w}}
\newcommand{\vx}{{\VectorFont x}}
\newcommand{\vy}{{\VectorFont y}}
\newcommand{\vz}{{\VectorFont z}}

\def\diag{\qopname\relax o{diag}}

\newtheorem{theorem}{Theorem}[section]
\newtheorem{lemma}[theorem]{Lemma}

\theoremstyle{definition}

\theoremstyle{remark}
\newtheorem{remark}[theorem]{Remark}

\theoremstyle{algorithm}
\newtheorem{algorithm}[theorem]{Algorithm}

\theoremstyle{corollary}

\theoremstyle{example}
\newtheorem{example}[theorem]{Example}

\theoremstyle{proposition}
\newtheorem{proposition}[theorem]{Proposition}

\usepackage{amssymb}
\usepackage{epsfig}
\usepackage{times}
\title{Dual Set Membership Filter with Minimizing  Nonlinear Transformation of Ellipsoid}

\author{Zhiguo Wang, Xiaojing Shen, Haiqi Liu, Fanqin Meng and Yunmin Zhu
\thanks{This work was supported  in part by   the special funds of NEDD of China under Grant No. 201314, the NSFC No. 61673282  and the PCSIRT16R53.}
\thanks{Zhiguo Wang is with School of Science and Engineering, The Chinese University of Hong Kong, Shenzhen, Guangdong, 518172, P.R. China, also with  Department of Electronic Engineering and Information Science, University of Science and Technology of China. Zhiguo Wang was with Department of Mathematics, Sichuan
University, Chengdu, Sichuan 610064, China.  Xiaojing Shen (corresponding author), Haiqi Liu, Fanqin Meng and Yunmin Zhu are with Department of Mathematics, Sichuan
University, Chengdu, Sichuan 610064, China. E-mail: wangzhiguo@cuhk.edu.cn, shenxj@scu.edu.cn, haiqiliu0330@163.com, mengfanqin2008@163.com, ymzhu@scu.edu.cn.}}

\begin{document}

 \maketitle
\begin{abstract}
In this paper, we propose a dual set membership filter for nonlinear dynamic systems with unknown but bounded noises, and it has three distinctive properties. Firstly,  the nonlinear system is translated into the linear system by leveraging a semi-infinite programming, rather than linearizing the nonlinear function. In fact, the semi-infinite programming is to find an  ellipsoid bounding the nonlinear transformation of an ellipsoid, which aims to compute a tight ellipsoid to cover the state. Secondly,  the  duality result of the semi-infinite programming is derived by a rigorous analysis, then a first order Frank-Wolfe method is developed to  efficiently solve it with a lower computation complexity.  Thirdly, the proposed filter can take advantage of the linear set membership filter framework and can work on-line without solving the semidefinite programming problem. Furthermore, we apply the dual set membership filter to a  typical scenario of mobile robot localization. Finally, two illustrative examples in the simulations show the advantages and effectiveness of the dual set membership filter.
\end{abstract}

\noindent{\bf keywords:}  Set membership filter, nonlinear system, unknown but bounded noise, Frank-Wolfe method, mobile robot localization.

\section{Introduction}\label{sec_1}
The nonlinear filter is an important research problem in
many fields, such as target tracking \cite{Lan-Li17}, navigation \cite{BarShalom-Li-Kirubarajan01} and mobile robot \cite{Hugh-Tim06}, etc. It is extensively researched in Bayesian framework, which is  based on stochastic assumptions about the process and measurement noises. If we know the exact distribution about these noises, then the classic nonlinear Bayesian filter, such as, extended Kalman filter (EKF) \cite{Masazade-Fardad-Varshney12}, unscented Kalman filter (UKF) \cite{Julier-04}, particle filter (PF) \cite{Arulampalam-Maskell-Gordon-Clapp02}, can obtain the better estimation. However, the statistical properties of  the process noise and measurement noise may be imprecisely known, which can lead to  degrade the Bayesian filter performance. It then seems more natural to assume that the state perturbations
and measurement noise are unknown but bounded \cite{Polyak-Nazin-Durieu-Walter04}.

In this paper, we consider the problem of ellipsoidal set membership filter (SMF) for the nonlinear dynamic systems with unknown but bounded noises, which does not require any assumption on the noise statistics.  The set-membership filter for linear systems was firstly proposed by Schweppe in \cite{Schweppe68}, and its basic idea is  propagating bounding ellipsoids \cite{Kurzhanski-valyi97} (or polyhedron, boxes, zonotopes) for dynamic systems with bounded noises. Recently, SMF has also been extensively explored, see \cite{Wu-Wang-Ye13,Rohou-Jaulin-Mihaylova-Bars17,ElGhaoui-Calafiore01,Shen-Zhu-Song-Luo11,Durieu-Walter-Polyak01}
 and references therein.

It is difficult to extend SMF to nonlinear dynamic systems, especially for on-line implementation. The reason is that a general nonlinear function maps an ellipsoid to an irregular set, which makes the set operation very complicated. Some researchers have considered the set membership filter for nonlinear systems \cite{Wang-Shen-Zhu19,Leong-Nair16,Gu-He-Han15,Sun-Alkhatib-Kargoll-Kreinovich-Neumann18,Chen-Hu18}. In \cite{Scholte-Campbell03}, the authors develop the extended set membership filter (ESMF) for a
general class of nonlinear systems  with on-line usage. Specifically, the nonlinear dynamics are linearized about the current
estimate, then the higher order remainder terms are then bounded by interval mathematics \cite{Moore66}, and the
remainder bounds are incorporated as additions to the process or sensor noise bounds, finally the classis SMF for linear system can be used. In \cite{Yang-Li10}, the authors employ the fuzzy modeling approach to approximate the nonlinear systems, and S-procedure technique to determine a state estimation ellipsoid. In \cite{Calafiore05}, the author proposes a nonlinear set membership filter (NSMF) based on a two-step prediction-correction. Each step of this filter requires to solve a semidefinite optimization problem (SDP) to  obtain the optimal outer-bounding
ellipsoid.

The  limitation of current SMF approaches for nonlinear system is that it needs to linearize the nonlinear function. However, linearized transformations are only reliable if the nonlinear function can be well approximated by a
linear function \cite{Julier-04}. Otherwise, for a general nonlinear function, the linearized approximation can be extremely poor \cite{Liu-Li15}, which leads to a bigger ellipsoid to bound the higher order remainder. In fact, the linearization  is to find an ellipsoid to approximate an irregular set, which is the nonlinear transformation of an ellipsoid. If we can directly calculate a minimum size ellipsoid to cover the  irregular set, the performance of SMF can be better (see Fig. \ref{fig_01}). Thus, these facts motivate us to consider nonlinear set-membership filter without linearizing the nonlinear functions.

On the other hand, the current SMF is that each step needs to solve an SDP problem \cite{ElGhaoui-Calafiore01,Ma-Wang-Lam-Kyriakoulis17}, of course, it can be solved by an interior point method. However, as problems have grown in size, such as large scale problems in multiple target tracking, it does not work well due to the higher computation complexity \cite{Ahipasaoglu-Sun-Todd08}. Hence, people have renewed their interest in the first-order method, which can be well extended to such a large scale problem.  In addition, Frank-Wolfe (FW) method is a popular first order method, and it can efficiently solve the minimum-volume enclosing ellipsoid problem for the convex hull of a finite set of points \cite{Todd16}. In this paper, we consider developing FW method to SMF, which makes each step of SMF very cheap to perform.

The main contribution of the paper is as follows. Based on the advanced optimization technique, we propose a dual set membership filter (DSMF) to recursively compute a bounding ellipsoid  to cover the state, and it has three vital edges.
\begin{itemize}
  \item Comparing with the traditional nonlinear SMF,  we use the semi-infinite programming to translate the nonlinear system into the linear system, rather than linearizing the nonlinear function.
  \item The dual problem of the semi-infinite programming is derived by  a rigorous analysis, and it can be solved  by the Frank-Wolfe method  with $\mathcal{O}(n^{2}+(n+1)m)$ arithmetic operations, which is lower than the  SDP  method $\mathcal{O}(n^5m^{1.5})$.
  \item The proposed filter can take advantage of the linear set membership filter framework, and it can be used on-line without assuming the bounds of the remainder by  linearizing the nonlinear function.
\end{itemize}
 Furthermore, we apply the proposed filter to a  typical scenario of mobile robot localization. Finally, two illustrative examples in the simulations show that the proposed filter can obtain a tighter ellipsoid and better estimation.

The rest of this paper is organized as follows. The problem of SMF for nonlinear system is formulated in Section \ref{sec_2}. In Section \ref{sec_3}, a novel DSMF is developed for computing the state estimation ellipsoid. The semi-infinite programming in DSMF is solved efficiently by Frank-Wolfe method in Section \ref{sec_4}. In Section \ref{sec_5}, we show the application of the proposed DSMF. Simulations and conclusions are given in Section  \ref{sec_6} and Section  \ref{sec_7}, respectively. Some technical proofs are provided in the appendix.

\emph{Notations:} For a square matrix $\mX\succeq0$  (resp. $\mX\succ0$ ) denotes $\mX$ is semidefinite (resp. positive-definite). The superscript $``T"$ denotes the transpose. $``\partial"$ is the gradient operator. $\diag(\cdot)$ and $\mI$ stand for a block diagonal matrix and the identity matrix with an appropriate dimension, respectively. Ellipsoid is described as $\mathcal {E}=\{\vx\in \mathcal {R}^n: (\vx-\hat{\vx})^T\mP^{-1}(\vx-\hat{\vx})\leq1\}=
\{\vx: \vx=\hat{\vx}+\mE\vu,\mP=\mE\mE^T, \|\vu\|\leq 1\}$, where $\mP$  and  $\hat{\vx}$ are the shape matrix and the center of the ellipsoid $\mathcal {E}$. The ``size" of the ellipsoid is the function of the shape matrix $\mP$, and it is denoted by $f(\mP)$. In this paper, $f(\mP)$ is either $logdet(\mP)$, which corresponds to the volume of the ellipsoid $\mathcal {E}$, or $tr(\mP)$, which means the sum of squares of the semiaxes lengths of the ellipsoid $\mathcal {E}$.

\section{Problem Formulation}\label{sec_2}
In target tracking, one of the major objective is to estimate the state trajectory of a target. In general, the
dynamic model for target tracking  describes the evolution of the
target vector with respect to time. The most common state-space models are given as follows
\begin{align}
\label{equ_1}\vx_{k+1}&=f_k(\vx_k)+\vw_k,\\
\label{equ_2}\vy_k&=h_k(\vx_k)+\vv_k,
\end{align}
where $\vx_k\in \mathbb{R}^n$ and $\vy_k\in \mathbb{R}^l$ are the target state and observation, respectively. $f_k(\vx_k)$ and $h_k(\vx_k)$ are nonlinear functions of state $\vx_k$. Here, $\vw_k$ and $\vv_k$ are the process and observation noises, respectively.

The uncertain process noise $\vw_k$ and measurement noise $\vv_k$  are assumed to be bounded by the following ellipsoids
\begin{align}
\label{equ_w}\mathcal{W}_k&=\{\vw_k:\vw_k^T\mQ_k^{-1}\vw_k\leq1\},\\
\label{equ_v}\mathcal{V}_k&=\{\vv_k:\vv_k^T\mR_k^{-1}\vv_k\leq1\},
\end{align}
where $\mQ_k$ and $\mR_k$ are the shape matrix of the ellipsoids $\mathcal{W}_k$ and $\mathcal{V}_k$, respectively. Both of them are known symmetric positive-definite matrices.

Suppose that the initial state $\vx_0$ belongs to a given bounding ellipsoid:
\begin{align}
\label{equ_3} \mathcal {E}_0&=\{\vx\in
R^n:(\vx-\hat{\vx}_0)^T\mP_0^{-1}(\vx-\hat{\vx}_0)\leq1\},
\end{align}
where  $\hat{\vx}_0$ is the center of the ellipsoid $\mathcal {E}_0$; $\mP_0$ is the shape matrix of the ellipsoid $\mathcal {E}_0$ which is a known symmetric positive-definite matrix.

The goal of this paper is to  obtain an ellipsoid $\mathcal {E}_{k+1}$ to bound the state $\vx_{k+1}$ at time $k+1$ by a recursion method. Specifically, at time $k$, assume that the set membership filter has obtained an ellipsoid $\mathcal {E}_k$ contains the state $\vx_k$, i.e.,
\begin{align}
\nonumber\mathcal {E}_k&=\{\vx:(\vx-\hat{\vx}_k)^T\mP_k^{-1}(\vx-\hat{\vx}_k)\leq1\}\qquad\qquad\\
\label{equ_4}&=\{\vx: \vx=\hat{\vx}_k+\mE_k\eta_k, \mP_k=\mE_k\mE_k^T, \parallel\eta_k\parallel\leq1\},
\end{align}
where $\hat{\vx}_k$ is the center of the ellipsoid $\mathcal {E}_k$; $\mP_k$ is a known symmetric positive-definite matrix. Based on the ellipsoid $\mathcal {E}_k$ and nonlinear state function, we derive a predicted ellipsoid $\mathcal {E}_{k+1|k}$ in the prediction step, which is
\begin{align}
\nonumber \mathcal{E}_{k+1|k}&=\{\vx:(\vx-\hat{\vx}_{k+1|k})^T\mP_{k+1|k}^{-1}
(\vx-\hat{\vx}_{k+1|k})\leq1\}\\
\label{equ_51a}& =\{\vx: \vx=\hat{\vx}_{k+1|k}+\mE_{k+1|k}\eta_{k+1|k}, \mP_{k+1|k}=\mE_{k+1|k}\mE_{k+1|k}^T, \parallel\eta_{k+1|k}\parallel\leq1\},
\end{align}
where $\hat{\vx}_{k+1|k}$ is the center of the ellipsoid $\mathcal {E}_{k+1|k}$ and $\mP_{k+1|k}$ is a symmetric positive-definite matrix.
 Then we use the  nonlinear measurement function (\ref{equ_2}) and measurement $\vy_k$ to obtain the updated ellipsoid $\mathcal {E}_{k+1}$  at time $k+1$ in the measurement update step, which is defined as follows
\begin{align}
\nonumber\mathcal {E}_{k+1}&=\{\vx:(\vx-\hat{\vx}_{k+1})^T\mP_{k+1}^{-1}(\vx-\hat{\vx}_{k+1})\leq1\}\qquad\qquad\\
\label{equ_4e}&=\{\vx: \vx=\hat{\vx}_{k+1}+\mE_{k+1}\eta_{k+1}, \mP_{k+1}=\mE_{k+1}\mE_{k+1}^T, \parallel\eta_{k+1}\parallel\leq1\},
\end{align}
where $\hat{\vx}_{k+1}$ is the center of the ellipsoid $\mathcal {E}_{k+1}$ and $\mP_{k+1}$ is a symmetric positive-definite matrix.

To  obtain a tighter ellipsoid, in next section, we develop a novel set membership filter by using an advanced optimization method to deal with the nonlinear function, rather than linearizing it.

\section{dual set Membership Filter Without Linearization}\label{sec_3}

In this section, a dual set membership filter is proposed for nonlinear system, which does not need to linearize the nonlinear function, and it  derives a predicted ellipsoid and an updated ellipsoid by solving the semi-infinite programming, respectively.

\subsection{Prediction Step}
Now, we consider the prediction step in DSMF, specifically, based on the ellipsoid $\mathcal{E}_k$ and the state equation at time $k$, we determine a predicted ellipsoid $\mathcal{E}_{k+1|k}$ that covers the set of state at time $k+1$. In general, there exist many ellipsoids containing the reachable set of states, however, finding a tighter predicted ellipsoid is difficult, especially in the nonlinear system.

The traditional ESMF \cite{Scholte-Campbell03} is linearizing the function $f_k$ in (\ref{equ_1}) about the current state estimate $\hat{\vx}_k$ (defined in (\ref{equ_4})) yields
\begin{align}
\label{equ_11aa}\vx_{k+1}=f_k(\hat{\vx}_k)+\frac{\partial f(\vx_k)}{\partial \vx}\Big|_{\vx_k=\hat{\vx}_k}(\vx_k-\hat{\vx}_k)+R_k^f(\vx_k,\hat{\vx}_k)+\vw_k,
\end{align}
where $R_k^f(\vx_k,\hat{\vx}_k)$ is the high-order remainder, and is written as
\begin{align}
\label{equ_12aa}R_k^f(\vx_k,\hat{\vx}_k)=(\vx_k-\hat{\vx}_k)^T\frac{1}{2}\frac{\partial^2 f(\bar{X}_{k})}{\partial \vx}
(\vx_k-\hat{\vx}_k), \\
\nonumber~\forall~\bar{X}_{k}=\hat{\vx}_k+\theta_k(\vx_k-\hat{\vx}_k), ~0\leq\theta_k\leq 1.
\end{align}
The other form for the remainder \cite{Wang-Shen-Zhu-Pan18} is
\begin{align}
\label{equ_13aa}R_k^f(\vx_k,\hat{\vx}_k)=f_k(\vx_k)-f_k(\hat{\vx}_k)-\frac{\partial f(\vx_k)}{\partial \vx}\Big|_{\vx_k=\hat{\vx}_k}(\vx_k-\hat{\vx}_k).
\end{align}

In \cite{Scholte-Campbell03} ,  the higher-order term $R_k^f(\vx_k,\hat{\vx}_k)$ in (\ref{equ_12aa}) can be bounded by an interval, then the interval is bounded by an ellipsoid. But this method is  conservative, and sometimes it is difficult  to calculate the Hessian matrix. For the second form in  (\ref{equ_13aa}), it is hard to analyze  its property and find a tighter ellipsoid to contain the remainder.

\begin{figure}[ht]
\vbox to 3cm{\vfill \hbox to \hsize{\hfill
\scalebox{0.5}[0.5]{\includegraphics{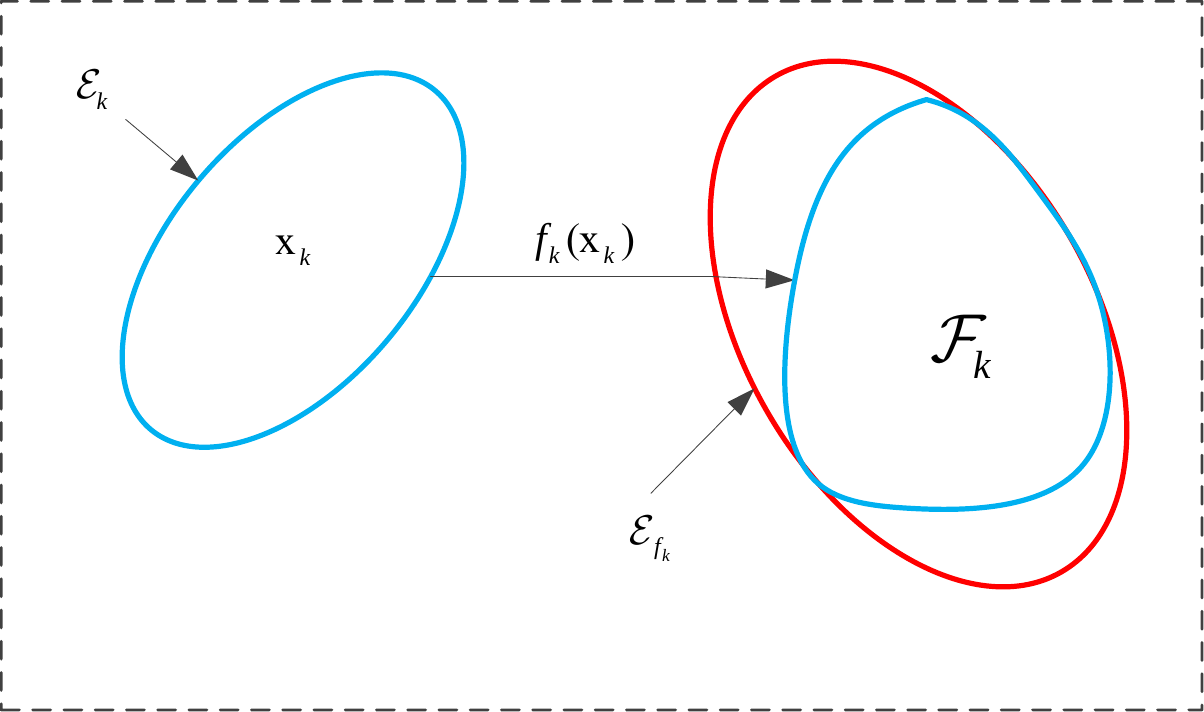}} \hfill}\vfill}
\caption{An illustration of the optimization problem (\ref{equ_49})-(\ref{equ_50}).}\label{fig_01}
\end{figure}

In this section, we use an optimization technique to avoid the above difficulty. Note that if the nonlinear state transform function $f_k$ is continuous, then $\mathcal{F}_k=\{f_k(\vx_k):\vx_k\in\mathcal{E}_k\}$ is a compact set. Comparing with linearizing the nonlinear function $f_k$, we directly derive  an ellipsoid $\mathcal{E}_{f_k}$ containing the nonlinear transformation  $\mathcal{F}_k$, i.e., $\mathcal{E}_{f_k}\supseteq \mathcal{F}_k$, by solving the following optimization problem
\begin{eqnarray}
\label{equ_49} &&\min ~f\big(\mP_{f_k}\big)\\
\label{equ_50}  && \mbox{s.t.}~ \big[\vx_{f_k}-\hat{\vx}_{f_k}\big]^T\mP_{f_k}^{-1}\big[\vx_{f_k}-\hat{\vx}_{f_k}\big]\leq 1,
 ~\forall~ \vx_{f_k}\in\mathcal{F}_k,
\end{eqnarray}
with variables $\hat{\vx}_{f_k}$ and $\mP_{f_k}$, where $\hat{\vx}_{f_k}$ and $\mP_{f_k}$ are the center and shape matrix of the ellipsoid $\mathcal{E}_{f_k}$, respectively. Fig. \ref{fig_01} gives us an illustration of the optimization problem (\ref{equ_49})-(\ref{equ_50}).  The optimization problem has the property that there are two variables appearing in infinitely many constraints, then it is called the semi-infinite programming. In the next section, we will provide some methods to solve the the optimization problem (\ref{equ_49})-(\ref{equ_50}) and  obtain the ellipsoid $\mathcal{E}_{f_k}$. Thus, we assume that we have obtained the ellipsoid $\mathcal{E}_{f_k}$ in this subsection.

As well known, the another form of the ellipsoid $\mathcal{E}_{f_k}$ is
\begin{eqnarray}
\nonumber \mathcal{E}_{f_k}=\{\vx: \vx=\hat{\vx}_{f_k}+\mE_{f_k}\eta_{f_k}, \|\eta_{f_k}\|\leq 1, \mE_{f_k}\mE_{f_k}^T=\mP_{f_k}\},
\end{eqnarray}
then there exists $\eta_{f_k}$ such that $f_k(\vx_k)=\hat{\vx}_{f_k}+\mE_{f_k}\eta_{f_k}$ by $f_k(\vx_k)\in \mathcal{E}_{f_k}$. Therefore, the nonlinear equation (\ref{equ_1}) can be translated to a linear equation as follows
\begin{eqnarray}
\label{equ_56}\vx_{k+1}=\hat{\vx}_{f_k}+\mE_{f_k}\eta_{f_k}+\vw_k.
\end{eqnarray}
According to equation (\ref{equ_56}), we assert that
\begin{eqnarray}
\label{equ_51} \vx_{k+1}\in\mathcal{E}_{f_k}\oplus\mathcal{W}_k.
\end{eqnarray}
 Although the  Minkowski sum $\mathcal{E}_{f_k}\oplus\mathcal{W}_k$ is not an ellipsoid,  in \cite{Kurzhanski-valyi97}, the authors have shown that there exists the  ellipsoid $\mathcal{E}_{k+1|k}$, that is an external approximation of the Minkowski sum, i.e,
\begin{eqnarray}
\label{equ_52} \mathcal{E}_{f_k}\oplus\mathcal{W}_k\subseteq\mathcal{E}_{k+1|k},
\end{eqnarray}
where the center and shape matrix of the predicted ellipsoid $\mathcal{E}_{k+1|k}$ are calculated as follows
\begin{eqnarray}
\label{equ_57} \hat{\vx}_{k+1|k}&=&\hat{\vx}_{f_k}\\
\label{equ_53} \mP_{k+1|k}(p)&=&(1+p^{-1})\mP_{f_k}+(1+p)\mQ_k
\end{eqnarray}
for  any $p>0$, respectively. Now, we may select an optimal external ellipsoid $\mathcal{E}_{k+1|k}$ relative to some criteria by solving
\begin{eqnarray}
\label{equ_o1} \min_{p>0} ~f\big(\mP_{k+1|k}(p)\big).
\end{eqnarray}

Fortunately, there exits a unique ellipsoid with minimal sum of squares of semiaxes containing the Minkowski sum \cite{Durieu-Walter-Polyak01}, where the optimal value of $p$ is defined by $p^{*}$, and
\begin{eqnarray}
\label{equ_54} p^{*}=\frac{\sqrt{tr(\mP_{f_k})}}{\sqrt{tr(\mQ_k)}}.
\end{eqnarray}
Let $\mP_{k+1|k}$ denote as $\mP_{k+1|k}(p^{*})$, then we have found the optimal predicted ellipsoid $\mathcal{E}_{k+1|k}$.

\subsection{Measurement Update Step}
In the measurement update step, we want to determine a minimal updated ellipsoid ${\mathcal{E}}_{k+1}$ that contains the states, which is consistent
with both the prediction $\mathcal{E}_{k+1|k}$  and the measurement at time $k+1$.

If we linearize the nonlinear measurement function $h_{k+1}$ about the current state estimate $\hat{\vx}_{k+1|k}$, it yields
\begin{align}
\label{equ_haa}\vy_{k+1}=h_{k+1}(\hat{\vx}_{k+1|k})+\frac{\partial h_{k+1}(\vx_{k+1})}{\partial \vx}\Big|_{\vx_{k+1}=\hat{\vx}_{k+1|k}}(\vx_{k+1}-\hat{\vx}_{k+1|k})
+R_{k+1}^h(\vx_{k+1},\hat{\vx}_{k+1|k})+\vw_{k+1},
\end{align}
where $R_{k+1}^h(\vx_{k+1},\hat{\vx}_{k+1|k})$ is the high-order remainder, and is written as
\begin{align}
\label{equ_haa2}R_{k+1}^h(\vx_{k+1},\hat{\vx}_{k+1|k})&=(\vx_{k+1}-\hat{\vx}_{k+1|k})^T\frac{1}{2}\frac{\partial^2 f(\bar{X}_{k+1})}{\partial \vx}
(\vx_{k+1}-\hat{\vx}_{k+1|k}), \\
\nonumber~\forall~\bar{X}_{k+1|k}&=\hat{\vx}_{k+1|k}+\theta_{k+1|k}(\vx_{k+1}-\hat{\vx}_{k+1|k}), ~0\leq\theta_{k+1|k}\leq 1.
\end{align}
According to (\ref{equ_51a}), the other form of the high-order remainder is
\begin{align}
\label{equ_haa3}R_{k+1}^h(\vx_{k+1},\hat{\vx}_{k+1|k})
&= \eta_{k+1|k}^T\mE_{k+1|k}^T\frac{1}{2}\frac{\partial^2 f(\bar{X}_{k+1|k})}{\partial \vx}\mE_{k+1|k}\eta_{k+1|k}.
\end{align}

Since the larger predicted ellipsoid means that the bigger shape matrix $\mE_{k+1|k}$, from (\ref{equ_haa3}), it shows that the remainder becomes larger with the larger predicted ellipsoid, which needs to find a larger measurement ellipsoid covering the remainder. Finally, the intersection of the  measurement ellipsoid and the predicted ellipsoid is larger, which  brings  a worse updated ellipsoid to bound the intersection. To obtain a better result, we can solve a semi-infinite optimization problem to cover the nonlinear transformation of predicted ellipsoid directly, rather than linearizing the nonlinear measurement function.

For the convenience of analysis, assume the nonlinear function $h_{k+1}$ exists a continuous inverse function $h_{k+1}$, and we can relax this assumption in Section \ref{sec_5}. Then the equation (\ref{equ_2}) can be rewritten as follows
\begin{eqnarray}%
\label{equ_12}h_{k+1}^{{-1}}(\vy_{k+1}-\vv_{k+1})=\mE_p\vx_{k+1},
\end{eqnarray}
where $\mE_p$ is the projection matrix. For example, when the measurement function \cite{Pu-Liu-Yan-Liu-Luo18} is denoted as
\begin{align}
h(\vx)=\left[
                       \begin{array}{c}
                         r \\
                        \theta \\
                       \end{array}
                     \right]=
\left[
                       \begin{array}{c}
                         \sqrt{(x-a)^2+(y-b)^2} \\
                         \arctan \frac{y-b}{x-a} \\
                       \end{array}
                     \right]
\end{align}
where $\vx=[x,\dot{x},y,\dot{y},\ddot{x},\ddot{y}]$, $(x,y)$ is the target position and $(a,b)$ is the sensor position. Then
\begin{align}
h^{-1}=\left[
             \begin{array}{c}
               x \\
               y \\
             \end{array}
           \right]=
\left[
                       \begin{array}{c}
                         r \cos(\theta)+a \\
                         r \sin(\theta)+b \\
                       \end{array}
                     \right]=\mE_p\vx,
\end{align}
where
\begin{align}
\nonumber \mE_p=\left[
                \begin{array}{cccccc}
                  1 & 0 & 0 & 0 & 0& 0 \\
                  0 & 0 & 1 & 0 & 0 &0 \\
                \end{array}
              \right].
\end{align}

Based on the uncertain set $\vv_{k+1}$ of the measurement noise, we hope to find a minimum measurement ellipsoid $\mathcal {E}(\hat{\vz}_{k+1},\mP_{z_{k+1}})$ to contain the left side of the equation (\ref{equ_12}), which can be described as follows
\begin{eqnarray}
\label{equ_13} &&\min ~f(\mP_{z_{k+1}})\\
\label{equ_14}  && \mbox{s.t.}~ \mathcal {C}_{k+1}\subseteq \mathcal {E}(\hat{\vz}_{k+1},\mP_{z_{k+1}}),
\end{eqnarray}
where $\mathcal {C}_{k+1}=\{h_{k+1}^{{-1}}(\vy_{k+1}-\vv_{k+1}):\vv_{k+1}\in\mathcal{V}_{k+1}\}$ and $f(\mP_{z_{k+1}})$ is the objective function, which is the ``size" of the shape matrix $\mP_{z_{k+1}}$.

The optimization problem (\ref{equ_13})-(\ref{equ_14}) is equivalent to
\begin{eqnarray}
\label{equ_15} &&\min ~f(\mP_{z_{k+1}})\\
\label{equ_16}  && \mbox{s.t.}~ (\vz_{k+1}-\hat{\vz}_{k+1})^T\mP_{z_{k+1}}^{{-1}}(\vz_{k+1}-\hat{\vz}_{k+1})\leq 1, ~\forall~ \vz_{k+1}\in\mathcal {C}_{k+1},
\end{eqnarray}
which is also a semi-infinite optimization problem, and it has the similar form with the optimization problem (\ref{equ_49})-(\ref{equ_50}). Then we will solve it at next section.

Since $\mE_p\vx_{k+1}$  belongs to $\mathcal {E}(\hat{\vz}_{k+1},\mP_{z_{k+1}})$ by equation (\ref{equ_12}), and it is equivalent to
\begin{eqnarray}
\label{equ_34}-\hat{\vz}_{k+1}=-\mE_p\vx_{k+1}+\mE_{z_{k+1}}\eta_{k+1}, \parallel\eta_{k+1}\parallel\leq1,
\end{eqnarray}
where $\mP_{z_{k+1}}=\mE_{z_{k+1}}\mE_{z_{k+1}}^T$. The nonlinear measurement equation (\ref{equ_2}) can be approximated by the above linear equation.

The measurement update step is to find a minimum volume or trace ellipsoid  containing the intersection of predicted ellipsoid $\mathcal{E}(\hat{\vx}_{k+1|k},\mP_{k+1|k})$ and measurement ellipsoid $\mathcal {E}(\hat{\vz}_{k+1},\mP_{z_{k+1}})$. Based on the famous result for linear set membership filter in \cite{Schweppe68,Wang-Shen-Zhu18}, the formula of the measurement update step can be written as follows
\begin{align}
\label{equ_35}\hat{\vx}_{k+1}&=\hat{\vx}_{k+1|k}+\frac{\mP_{k+1|k}}{1-\rho_{k+1}}(\mE_{k+1}^p)^T\Big[\mE_{k+1}^p\frac{\mP_{k+1|k}}
{1-\rho_{k+1}}(\mE_{k+1}^p)^T+\frac{\mP_{z_{k+1}}}{\rho_{k+1}}\Big]^{-1}(\hat{\vz}_{k+1}-\mE_{k+1}^p\hat{\vx}_{k+1|k})\\
\label{equ_36} \bar{\mP}_{k+1}&=\big[(1-\rho_{k+1})\mP_{k+1|k}^{-1}+\rho_{k+1}(\mE_{k+1}^p)^T\mP_{z_{k+1}}^{-1}\mE_{k+1}^p\big]^{-1}, ~\rho_{k+1}\geq 0,\\
\label{equ_37} \delta_{k+1}&=(\hat{\vz}_{k+1}-\mE_{k+1}^p\hat{\vx}_{k+1|k})^T
\Big[\mE_{k+1}^p\frac{\mP_{k+1|k}}{1-\rho_{k+1}}(\mE_{k+1}^p)^T+\frac{\mP_{z_{k+1}}}{\rho_{k+1}}\Big]^{-1}(\hat{\vz}_{k+1}-\mE_{k+1}^p\hat{\vx}_{k+1|k})\\
\label{equ_38} \mP_{k+1}&=(1-\delta_{k+1})\bar{\mP}_{k+1}.
\end{align}

The optimal parameter $\rho_{k+1}^{*}$ can be obtained by solving the following problem
 \begin{eqnarray}
\label{equ_58} &&\min ~f(\mP_{k+1})\\
\label{equ_59}  && \mbox{s.t.}~ 1\geq\rho_{k+1}\geq 0.
\end{eqnarray}
Substituting $\rho_{k+1}^{*}$ to the equations (\ref{equ_35})-(\ref{equ_38}), we obtain the final center and shape matrix of the updated ellipsoid.
\begin{remark}
Searching for an outer-bounding ellipsoid ${\mathcal{E}}_{k+1}$ is reduced to a one-dimensional optimization by (\ref{equ_58})-(\ref{equ_59}). Meanwhile, the optimization problem (\ref{equ_58})-(\ref{equ_59}) is a convex optimization problem, which can be solved efficiently by many methods, such as the golden section method. In \cite{Calafiore05}, the author obtains the predicted ellipsoid and the updated ellipsoid by solving an SDP problem with $\mathcal{O}(n^3)$. Comparing with it, our method in each step has smaller calculation efforts with $\mathcal{O}(n^{2}+(n+1)m)$.
\end{remark}
\begin{example}
Let the sensor position be origin point. The center and  shape matrix  of the predicted ellipsoid are $[10~20]^T$  and $\sigma\mI$, respectively. Here, $\sigma$ is a parameter. The shape matrix of the measurement noise is $R=\diag(10,1)$. Assume the true state is sampling from the predicted ellipsoid, then we use the measurement equation to obtain a measurement $y$. According to the predicted ellipsoid and measurement, we exploit ESMF and the new method to get the updated ellipsoid $\mathcal{E}$, and its shape matrix is $\mP$. Figs. \ref{fig_02}-\ref{fig_03} show some results. Obviously, the new method can obtain a smaller updated ellipsoid in Fig. \ref{fig_02}. From Fig. \ref{fig_03}, we can see that with the increasing of $\sigma$, which means  the predicted ellipsoid  becomes larger, then the size of the estimated ellipsoid is raising. However, the new method gets a stable estimated ellipsoid. The reason may be that our method does not need to linearize the nonlinear measurement function, which does not bring much uncertainty.
\begin{figure}[ht]
\vbox to 5cm{\vfill \hbox to \hsize{\hfill
\scalebox{0.5}[0.5]{\includegraphics{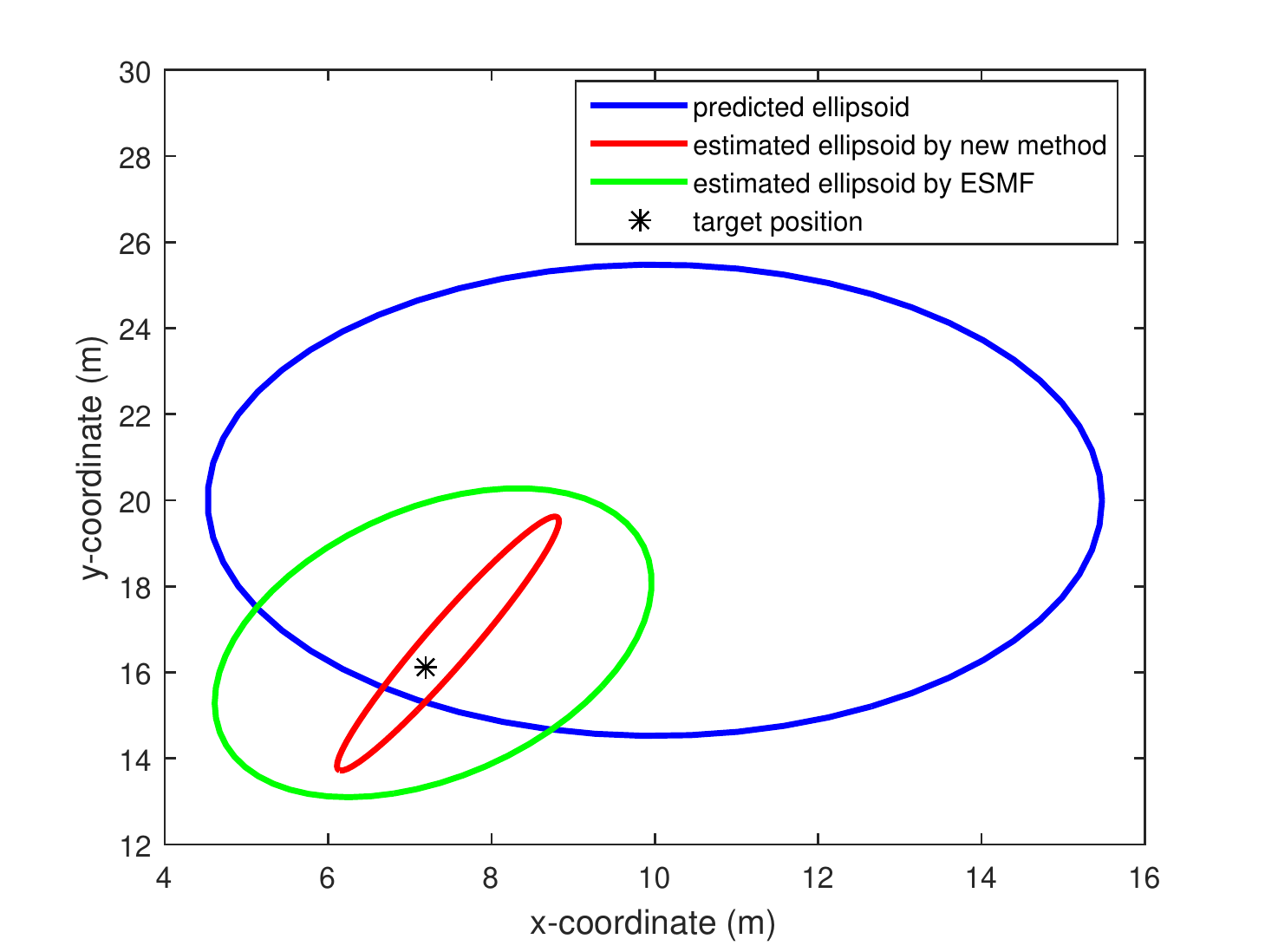}} \hfill}\vfill}
\caption{The estimated ellipsoid with the proposed method and ESMF with $\sigma=30$.}\label{fig_02}
\end{figure}
\begin{figure}[ht]
\vbox to 5cm{\vfill \hbox to \hsize{\hfill
\scalebox{0.5}[0.5]{\includegraphics{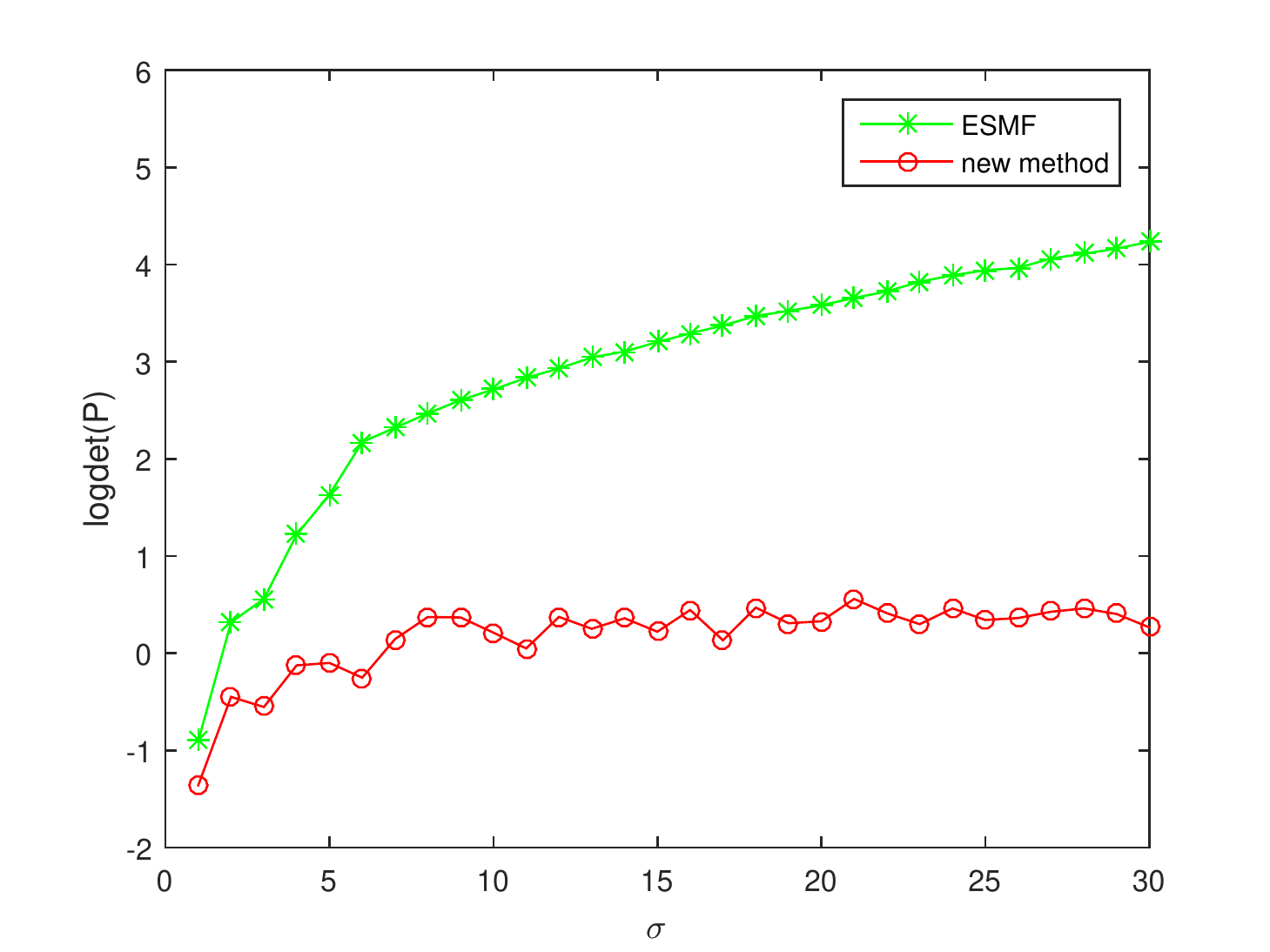}} \hfill}\vfill}
\caption{The size (logdet) of the estimated ellipsoid is plotted as a function of $\sigma$ for 50 Monte Carlo runs.}\label{fig_03}
\end{figure}

\end{example}
A summary of the proposed dual set membership filter is in Algorithm 1, and we can apply it to target tracking. In the next section, we provide an efficient method to solve the semi-infinite optimization problem in Algorithm 1.

\begin{algorithm}\label{alg_1}
[Dual Set Membership Filter]
\begin{itemize}
  \item \textbf{Input:} The nonlinear state function $f_k$, measurement function $h_k$,
the measurement $\vy_k$,
the center $\vx_0$ and shape matrix $\mP_0$ of initial ellipsoid;
  \item \textbf{Output:}
  \item 1.  For each $k=1:K$ do
  \item ~~~~2.  Solve optimization problem (\ref{equ_49})-(\ref{equ_50}) to obtain $\hat{\vx}_{f_k}$ and $\mP_{f_k}$ (see Algorithm 3);
  \item ~~~~3.  Calculate the predicted ellipsoid $\mathcal{E}_{k+1|k}$ by (\ref{equ_57})-(\ref{equ_53});
  \item ~~~~4.  Solve optimization problem (\ref{equ_15})-(\ref{equ_16}) to obtain $\hat{\vz}_{k}$ and $\mP_{z_k}$ (see Algorithm 3);
  \item ~~~~5.  Calculate the updated ellipsoid $\mathcal{E}_{k+1}$ by (\ref{equ_35})-(\ref{equ_38});
  \item 6.  \textbf{End For} 
  \item \textbf{Return} $\hat{\vx}_{k+1}$ and $\mP_{k+1}$;
\end{itemize}
\end{algorithm}

\section{Solving  the semi-infinite programming}\label{sec_4}
To overcome the difficult of the nonlinear function, we have translated the nonlinear dynamic system into the linear system (\ref{equ_56}) and (\ref{equ_34}) by solving the semi-infinite programming (\ref{equ_49})-(\ref{equ_50}) and (\ref{equ_15})-(\ref{equ_16}). In this section, we focus on discussing how to solve the semi-infinite programming.

Firstly, the optimization problems (\ref{equ_49})-(\ref{equ_50}) and (\ref{equ_15})-(\ref{equ_16}) can be unified as
\begin{eqnarray}
\label{equ_61} &&\min_{\mP,\hat{\vy}} ~logdet(\mP)\\
\label{equ_62}  && \mbox{s.t.}~ (\vy-\hat{\vy})^T\mP^{-1}(\vy-\hat{\vy})\leq n,
 ~\forall~ \vy\in\mathcal{C}.
\end{eqnarray}
Here the vector $\hat{\vy}$ and matrix $\mP$ are the optimization variables, and they are also the center and the shape matrix of ellipsoid $\mathcal{E}$. $n$ is the dimension of $\vy$. The set $\mathcal{C}$  is defined as
\begin{eqnarray}
\label{equ_62a}\mathcal{C}&=&\{\vy:\vy=f(\vx),~\vx\in\mathcal{E}_0\},\\
\label{equ_62aa}\mathcal{E}_0&=&\{\vy: (\vy-\hat{\vy}_0)^T\mP_0^{-1}(\vy-\hat{\vy}_0)\leq 1\},
\end{eqnarray}
where $f$ is a nonlinear continuous function, $\hat{\vy}_0$ and $\mP_0$ are the center and shape matrix of the ellipsoid $\mathcal{E}_0$, respectively. If the objective function (\ref{equ_61}) is the trace function, i.e., $tr(\mP)$, the similar result and algorithm can be obtained.

\subsection{A First Order Optimal Condition}
In this subsection, we discuss the problem (\ref{equ_61})-(\ref{equ_62}) from a theoretical point of view, and its dual problem helps us to design a fast algorithm to solve the problem (\ref{equ_61})-(\ref{equ_62}). Thus, we call the proposed filter as dual set membership filter.  The following proposition shows that strong duality holds.
\begin{proposition}\label{pro_1}
The optimal solutions of problem (\ref{equ_61})-(\ref{equ_62}) are
\begin{eqnarray}
\hat{\vy}^{*}=\int_{\mathcal{C}}\vy d\mu^{*},~ \mP^{*}=\int_{\mathcal{C}}\vy \vy^T d\mu^{*}-\hat{\vy}^{*}\hat{\vy}^{{*}^T},
\end{eqnarray}
where $\mu^{*}$ is a measure and the optimal solution of the following optimization problem
\begin{eqnarray}
 \label{equ_84}&&\max_{\mu} ~\mbox{logdet}\left(\int_{\mathcal{C}\times \{1\}}\tilde{\vy}\tilde{\vy}^T d\mu\right)\\
 \label{equ_85}&& \mbox{s.t.}~~ \int_{\mathcal{C}\times \{1\}}d\mu=1, \mu\geq 0,
\end{eqnarray}
where $\tilde{\vy}=[\vy^T,1]^T$.
\end{proposition}
\begin{proof}
See Appendix.
\end{proof}
\begin{remark}
Comparing the optimization problems (\ref{equ_61})-(\ref{equ_62}) with (\ref{equ_84})-(\ref{equ_85}), the constraint function in (\ref{equ_85}) is a linear function on a continuous measure $\mu$, which is simpler to find the feasible solution.  If $\mu$ is a discrete measure, the constrain (\ref{equ_85}) is a simplex, which offers guidance on the efficient first order algorithm.
\end{remark}
Since the proof of Proposition \ref{pro_1} needs the following two lemmas, and these two lemmas are important to our algorithm, we firstly show them as follows.

\begin{lemma}\label{lem_1} Let $d=n+1$, the minimum of problem
\begin{eqnarray}
\label{equ_63}&&\min_{\mH} ~-\mbox{logdet}(\mH)\\
\label{equ_64}&& \mbox{s.t.}~ \tilde{\vy}^T\mH\tilde{\vy}\leq d, ~\forall ~\tilde{\vy}\in \Omega={\mathcal{C}\times \{1\}},
\end{eqnarray}
is obtained at $\mH^{*}$ if and only if $\mH^{*}$  is feasible and there exists a $\mu^{*}$ such that
\begin{eqnarray}
\label{equ_69} {\mH^{*}}^{-1}=\int_{\Omega}\tilde{\vy}\tilde{\vy}^T d\mu^{*},\\
\label{equ_70} \int_{\Omega}(\tilde{\vy}^T\mH^{*}\tilde{\vy}-d)d\mu^{*}=0,\\
\label{equ_71} \mu^{*}\geq 0.
\end{eqnarray}
\end{lemma}
\begin{proof}
See Appendix.
\end{proof}
\begin{lemma}\label{lem_2}
The dual problem corresponding to the primal problem (\ref{equ_63})-(\ref{equ_64}) is
\begin{eqnarray}
\label{equ_72}&&\max_{\mu} ~\mbox{logdet}\left(\int_{\Omega}\tilde{\vy}\tilde{\vy}^T d\mu\right)\\
\label{equ_73}&& \mbox{s.t.}~ \int_{\Omega}\mu=1, \mu\geq 0.
\end{eqnarray}
\end{lemma}
\begin{proof}
See Appendix.
\end{proof}

In fact, it suffices to consider the discrete measures $\mu$ from John's optimality conditions \cite{Todd16}, which puts positive measure on a finite points of the set $\Omega$.  Then we randomly sample some points $\tilde{\vy}_i$ from the set $\Omega$ so that the optimization problem (\ref{equ_84})-(\ref{equ_85}) can be approximated by
\begin{eqnarray}
 \label{equ_86}&&\max_{\mu_i} ~\mbox{logdet}\left(\sum_{i=1}^m\mu_i\tilde{\vy}_i\tilde{\vy}_i^T\right)\\
 \label{equ_87}&& \mbox{s.t.}~~ \sum_{i=1}^m\mu_i=1, \mu_i\geq 0.
\end{eqnarray}
Since $\Omega=\mathcal{C}\times \{1\}$ and $\tilde{\vy}=[\vy^T,1]^T$, if we want to get points $\tilde{\vy}_i$ from the set $\Omega$, it is sufficient to get points $\vy_i$ from the set $\mathcal{C}$, and $i=1,\ldots,m$.
\subsection{The Relationship Between DSMF and UKF}
In this subsection, we firstly consider three cases for obtaining the points  $\vy_i$ from the set $\mathcal{C}$.
\begin{itemize}
  \item For a general function $f$ defined in (\ref{equ_62a}),  it is easy to sample the interior points or boundary points $\vx_i$ from the ellipsoid $\mathcal{E}_0$, then  the points in the set $\mathcal{C}$ can be obtained directly by $\vy_i=f(\vx_i)$, and $i=1,\ldots,m$.
  \item For some particular nonlinear functions, such as homeomorphic function \cite{Armstrong83},  it takes the boundary of the ellipsoid $\mathcal{E}_0$ to the boundary of the set $\mathcal{C}$. Since  $\int_{\Omega}(\tilde{\vy}^T\mH^{*}\tilde{\vy}-d)d\mu^{*}=0$ in (\ref{equ_70}), when $\tilde{\vy}^T\mH^{*}\tilde{\vy}=d$, $\mu\neq0$, otherwise $\mu=0$. In order to obtain the points that making $\mu\neq0$, it is sufficient for getting the boundary points of the set $\mathcal{C}$.
  \item For some common nonlinear functions in two or three dimensional radar system \cite{Wang-Shen-Zhu-Pan18}, we can  sample the boundary points  from the ellipsoid $\mathcal{E}_0$, then it is enough to use them to get the boundary points $\vy_i$ in $\mathcal{C}$ and $i=1,\ldots,m$.
\end{itemize}

Next, we show the relationship between DSMF and UKF. As we know, for UKF, a set of points (sigma points) are chosen and the nonlinear function is applied to these points, in
turn, to yield a cloud of transformed points. Then the statistics of
the transformed points can be approximated to form an estimate
of the nonlinearly transformed mean and covariance \cite{Julier-04}.

For our methods, a set of interior points or boundary points are chosen from the ellipsoid $\mathcal{E}_0$ and  the nonlinear function is also applied to them, however, we want to find an ellipsoid to cover all these nonlinear transformed points, which is formulated as the optimization problem (\ref{equ_86})-(\ref{equ_87}). Thus, comparing with UKF, DSMF can obtain a robust estimate of the true state  by solving the optimization problem (\ref{equ_86})-(\ref{equ_87}).

\subsection{Frank-Wolfe Method}
The Frank-Wolfe (FW) or conditional gradient algorithm is a projection-free  method, it solves the problem (\ref{equ_86})-(\ref{equ_87}) without requiring  to compute a projection onto the feasible set.

Let $g(\mu)=\mbox{logdet}(\sum_{i=1}^m\mu_i\tilde{\vy}_i\tilde{\vy}_i^T)$ and $M=\{\mu: \sum_{i=1}^m\mu_i=1, \mu_i\geq 0\}$. If we make a first-order linear (Taylor) approximation to the objective function $g(\mu)$ at the current solution $\mu^t$:
\begin{align}
g(\mu)\approx g(\mu^t)+\omega^{t^T}(\mu-\mu^t),
\end{align}
where $\omega^t:=\omega(\mu^t)=\partial g(\mu^t)$. We might then consider maximizing this linear function over the unit
simplex $M$, and  the analysis solution is at a vertex, i.e., unit
vector $e_i$, which means the $i$th component is 1, the other is 0. Hence we might
wish to move along the line from  $\mu^t$ towards $e_i$ to update $\mu^{t+1}$. The pseudo code for the Frank-Wolfe algorithm is shown in Algorithm 2.
\begin{algorithm}\label{alg_2}
[FW algorithm]
\begin{itemize}
  \item \textbf{Input: }Let $\mu^0\in M$
  \item \textbf{Output: }
  \item 1: For {each $t=1:T$}
  \item ~~~~2: Compute $e_i^{t}:\arg\min_{s\in M}\langle s,~ \partial g(\mu^t) \rangle$;
  \item ~~~~3: Let $d_t=e_i^{t}-\mu^t$;
  \item ~~~~4: Obtain the optimal step size $\gamma_t=\mathop{\arg\min}\limits_{\gamma\in[0,1]}f(\mu^t+\gamma d_t)$;
  \item ~~~~5: Update $\mu^{t+1}=\mu^t+\gamma_t d_t$.
  \item 6: \textbf{End For}
  \item \textbf{Return} $\mu^T$;
\end{itemize}
\end{algorithm}

In \cite{Ahipasaoglu-Sun-Todd08}, the authors have showed the linear convergence of the Frank-Wolfe algorithm for the optimization problem (\ref{equ_86})-(\ref{equ_87}), and it has a decisive advantage in a large-scale problem. Based on Proposition \ref{pro_1}, we can use the following Algorithm 3 to solve the semi-infinite programming (\ref{equ_49})-(\ref{equ_50}) and (\ref{equ_15})-(\ref{equ_16}).

\begin{algorithm}[htb]\label{alg_3}
[A First Order Method for Semi-Infinite Programming.]
\begin{itemize}
  
  \item \textbf{Input:} The number of samples $m$;  The set $\mathcal{C}$;
  \item \textbf{Output:} Obtain ellipsoid $\mathcal{E}$ to contain set $\mathcal{C}$;
  \item 1: Generate samples $\vy_1,\ldots,\vy_m$ from set $\mathcal{C}$, and let $\tilde{\vy}_i=[\vy_i^T,~1]^T$;
  \item 2: Solve the optimization problem (\ref{equ_86})-(\ref{equ_87}) to obtain $\mu_i^{*}$ by FW algorithm;
  \item 3: \textbf{Return:} $\hat{\vy}=\sum_{i=1}^m\mu_i^{*}\vy_i ,~ \mP=\sum_{i=1}^m\mu_i^{*}\vy_i \vy_i^T-\hat{\vy}\hat{\vy}^T$;
\end{itemize}
\end{algorithm}
\subsection{Frank-Wolfe Method versus SDP}
According to Schur complement, the original optimization problem (\ref{equ_61})-(\ref{equ_62}) by introducing variables $\mP=\mB^{-1}$ and $\hat{\vy}=\mP\hat{\vb}$ can be equivalent to
\begin{eqnarray}
\label{equ_88} &&\min_{\mB,\hat{\vb}} ~-\mbox{logdet}(\mB)\\
\label{equ_89}  && \mbox{s.t.}~ \left[
                                  \begin{array}{cc}
                                    n & \vy^T\mB-\hat{\vb}^T \\
                                    \mB\vy-\hat{\vb} & \mB \\
                                  \end{array}
                                \right]\geq 0,
 ~\forall~ \vy\in\mathcal{C},
\end{eqnarray}
where  $\mB$ and $\hat{\vb}$ are the optimal variable. In order to solve the problem (\ref{equ_88})-(\ref{equ_89}), we can take samples from the boundary and the interior of the set  $\mathcal{C}$, so that we can obtain a finite set of  $\vy_1,\ldots,\vy_m$, then the infinite constraint (\ref{equ_89}) can be approximated by $m$ constraints based on $\vy_1,\ldots,\vy_m$. Specifically, the optimization problem (\ref{equ_88})-(\ref{equ_89}) can be approximated by an SDP problem
\begin{eqnarray}
\label{equ_90} &&\min_{\mB,\hat{\vb}} ~-\mbox{logdet}(\mB)\\
\label{equ_91}  && \mbox{s.t.}~ \left[
                                  \begin{array}{cc}
                                    n & \vy_i^T\mB-\hat{\vb}^T \\
                                    \mB\vy_i-\hat{\vb} & \mB \\
                                  \end{array}
                                \right]\geq 0,
 ~\forall~ i=1,\ldots,m.
\end{eqnarray}
In this problem the variable is the matrix $\mB$ and $\hat{\vb}$, so the dimension of the optimization
variable is $N=n(n+1)/2+n$.  \cite{Calafiore-Campi05} used statistical learning techniques to provide an explicit bound on the measure of the set of original constraints that are possibly violated by the randomized sampling, and they proved that this measure rapidly decreases to zero as $m$ is increasing. Therefore, the solution of the optimization problem (\ref{equ_90})-(\ref{equ_91}) can be made approximately feasible for the semi-infinite optimization problem (\ref{equ_88})-(\ref{equ_89}) by sampling a sufficient number of constraints.

When the inequalities (\ref{equ_91}) are combined into one
large inequalities, the dimension of the constrain matrix is $m(n+1)$. Vandenberghe and Boyd \cite{Boyd-ElGhaoui-Feron-Balakrishnan94} have developed a (primal-dual) interior-point method
that solves (\ref{equ_90})-(\ref{equ_91}) by exploiting the problem structure. They prove the worst-case estimate of $\mathcal{O}(N^{2.75}m^{1.5})$ arithmetic operations to solve the problem to a given accuracy. However, each iteration in FW method requires $\mathcal{O}(n^{2}+(n+1)m)$ arithmetic operations \cite{Todd16}. The cheapness of the iterations in FW method, together with their (relatively)
attractive convergence properties \cite{Ahipasaoglu-Sun-Todd08}, leads to that it is efficiency for large number $m$.
\begin{example}
There are $m$ random vector $\vy_i, i=1,\ldots,m$ that are generated by the standard uniform distribution. Assume the dimension of $\vy_i$ is $n$. We solve the SDP problem (\ref{equ_90})-(\ref{equ_91}) and the optimization problem (\ref{equ_86})-(\ref{equ_87}) by SDPT3 algorithm using the CVX platform and FW algorithm, respectively. The results presented in Table \ref{tab_1}  are mean solution times of 20 Monte Carlo runs   by the SDP solver in the third column  and by FW algorithm in the forth column. It is shown that the first-order FW  algorithm dominates the SDP method, sometimes it is more than 200 times
faster. Thanks to the faster FW algorithm, the proposed set membership filter can work on-line.
\end{example}
\begin{table}[!hbp]
\centering
\caption{Mean solution running time of SDPT3 and FW} \label{tab_1}
\begin{tabular}{|c|c|c|c|}
  \hline
  n & m & SDP & FW \\\hline
  2 & 50 &  0.7236& 0.0027 \\
  2 & 100 & 1.2269& 0.0037 \\
  2 & 200 & 2.2153 & 0.0076 \\
  2 & 400 & 4.1819 & 0.0048 \\
  2 & 600 & 6.2677 & 0.0115 \\
  2 & 800 & 8.4777 & 0.0141 \\
  2 & 1000 & 10.7986 & 0.0070 \\
  6 & 50 &  0.8008& 0.0108 \\
  6 & 100 & 1.3663& 0.0119\\
  6 & 200 & 2.5366 & 0.0185 \\
  6 & 400 & 5.0538 & 0.0260 \\
  6 & 600 & 7.7717 & 0.0305 \\
  6 & 800 & 10.7293& 0.0331 \\
 6 & 1000 & 13.9558& 0.0416 \\
  \hline
\end{tabular}
\end{table}
\section{Applications of dual set-Membership Filter}\label{sec_5}
Mobile robot location is an important area in artificial intelligence, and one of the key problem about it is the simultaneous localization and mapping (SLAM). In the SLAM problem, mobile robots need maps to locate themselves, and maps need robots to update themselves \cite{Jaulin09,Chen-Hu-Ho-Yu18,Yu-Zamora-Soria16}.

Consider a mobile robot moving through planar environments at time $k$, the state vector of robot is defined as $\vx_k=[p_k^x,p_k^y,\theta_k]$, where $(p_k^x,p_k^y)$ describes robot's position in X-Y plan, and $\theta_k$ is to define the angular orientation. The nonlinear motion model of the robot is given by
\begin{align}
\label{equ_5_9} \vx_{k+1}=f_k(\vx_k)+\vw_k,
\end{align}
where $f_k$ is a nonlinear state function.
Assume the location of the landmark is $(s^x,s^y)$, then the measurement model is \cite{Thrun02}
\begin{align}
\label{equ_5_10} \vy_k=\left[
                         \begin{array}{c}
                           r_k \\
                           \phi_k \\
                         \end{array}
                       \right]=\left[
                                 \begin{array}{c}
                                   \sqrt{(p_k^x-s^x)^2+(p_k^y-s^y)^2} \\
                                   \theta_k-\arctan(\frac{p_k^y-s^y}{p_k^x-s^x}) \\
                                 \end{array}
                               \right]+\vv_k
.
\end{align}
In fact, the nonlinear function in the right side of the equation (\ref{equ_5_10}) does not exist a continuous inverse function. However, after some transformations, we can get
\begin{align}
\label{equ_5_11} \mE\vx_k=\left[
                   \begin{array}{c}
                     p_k^x \\
                     p_k^y \\
                   \end{array}
                 \right]=g(\vy_k,\vv_k,\theta_k)\triangleq\left[
                           \begin{array}{c}
                             (\vy_k^1-\vv_k^1)\cos(\theta_k-\phi_k-\vv_k^2)+s^x \\
                             (\vy_k^1-\vv_k^1)\sin(\theta_k-\phi_k-\vv_k^2)+s^y \\
                           \end{array}
                         \right],
\end{align}
where $\vy_k=[\vy_k^1~\vy_k^2]$ and $\vv_k=[\vv_k^1~\vv_k^2]$. $\mE$ is the projection matrix with $\mE=\left[
                                                                                                              \begin{array}{ccc}
                                                                                                                1 & 0 & 0 \\
                                                                                                                0 & 1 & 0 \\
                                                                                                              \end{array}
                                                                                                            \right]
$.

Since the nonlinear function $g$ not only depends on $\vy_k$ and $\vv_k$, but also $\theta_k$, we cannot obtain the measurement ellipsoid if we only use the information of the measurement $\vy_k$ and noise ellipsoid $\mathcal{V}_k$ just as (\ref{equ_12}). To overcome this, we use the predicted set about $\theta_k$ such that $\theta_k\in \mathcal{C}^{\theta}_k=\{\theta: |\theta-\hat{\theta}_{k|k-1}|\leq \mP_{k|k-1}^{3,3}\}$, where $\hat{\theta}_{k|k-1}$ is the third component of the center $\hat{\vx}_{k+1|k}$ of the ellipsoid $\mathcal{E}_{k+1|k}$, and $\mP_{k|k-1}^{3,3}$ is the third row and the third column of the shape matrix ${\mP}_{k+1|k}$ of the predicted ellipsoid $\mathcal{E}_{k+1|k}$. Then, the measurement ellipsoid $\mathcal {E}(\hat{\vz}_{k},\mP_{z_{k}})$ can be obtained by solving the following optimization problem
\begin{eqnarray}
\label{equ_5_12} &&\min ~f(\mP_{z_{k}})\\
\label{equ_5_13}  && \mbox{s.t.}~ (\vz_{k}-\hat{\vz}_{k})^T\mP_{z_{k}}^{{-1}}(\vz_{k}-\hat{\vz}_{k})\leq 1, ~\forall~ \vz_{k}\in\mathcal {C}_{k},
\end{eqnarray}
where $\mathcal {C}_{k}=\{g(\vy_k,\vv_k,\theta_k):\vv_{k}\in\mathcal{V}_{k},\theta_k\in \mathcal{C}^{\theta}_k\}$, which is different with the definition in (\ref{equ_14}).

However, we can also get interior points or boundary points of the set $\mathcal {C}_{k}$, then the optimization problem (\ref{equ_5_12})-(\ref{equ_5_13}) can be solved efficiently by FW algorithm.  Thus, the proposed dual set membership filter can be used to obtain the mobile robot localization.

\section{Simulations}\label{sec_6}
In this section, we compare the performance among DSMF and UKF \cite{Julier-04}, ESMF \cite{Scholte-Campbell03}, NSMF \cite{Calafiore05} by two typical examples, and they show the advantages and effectiveness of the proposed DSMF.
\subsection{Nonlinear Measurement System}
\begin{figure}[ht]
\vbox to 5cm{\vfill \hbox to \hsize{\hfill
\scalebox{0.5}[0.5]{\includegraphics{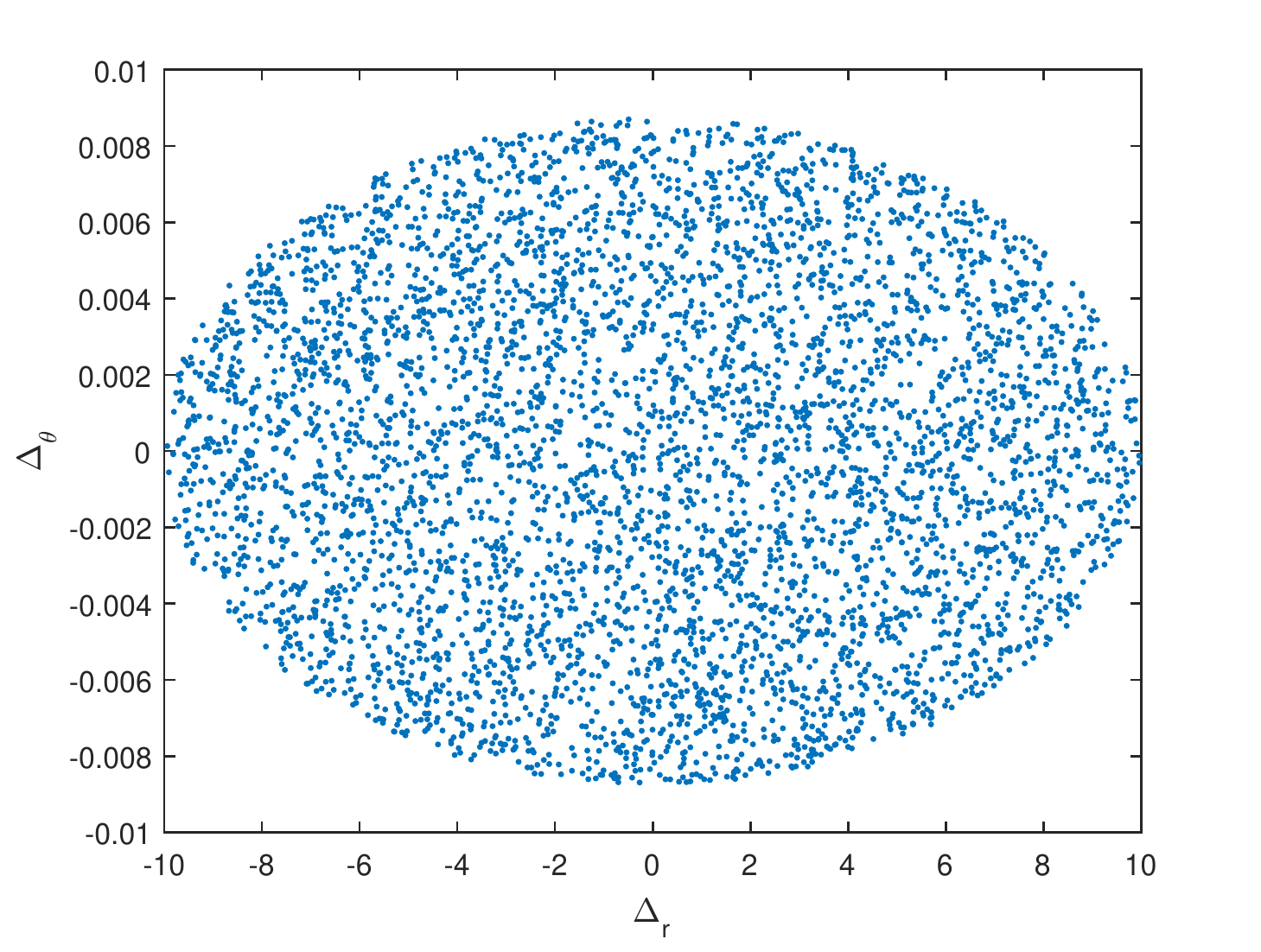}} \hfill}\vfill}
\caption{Measurement noise.}\label{fig_04}
\end{figure}

Consider the problem of
tracking a target in two dimensions. The state contains position and velocity
of $x$ and $y$ directions.
The dynamic system is
\begin{align}
\label{equ_46}\vx_{k+1}&=\mF_k\vx_k+\vw_k,\\
\label{equ_47}\vy_{k+1}&=h(\vx_{k+1})+\vv_{k+1}.
\end{align}
where
\begin{align}
\nonumber \mF_k= \left[
         \begin{array}{cccc}
           1 & 0 & T & 0 \\
           0 & 1 & 0 & T \\
           0 & 0 & 1 & 0 \\
           0 & 0 & 0 & 1 \\
         \end{array}
       \right],
                      h_{k}=
                      \left[
                       \begin{array}{c}
                         \sqrt{(p_k^x-420)^2+(p_k^y-420)^2} \\
                         \arctan \frac{p_k^y-420}{p_k^x-420} \\
                       \end{array}
                     \right].
\end{align}
Here, $T$ is the time sampling interval with $T=1$, $\vx_k$ is the state at time $k$, and $\vx_k=[p_k^x,p_k^y,v_k^x,v_k^y]$. The initial state is $\vx_0=[50 ~30~ 5~ 5]$.
Moreover, the process noise $\vw_k$ and measurement noise $\vv_k$ are taking value in specified ellipsoidal sets $\mathcal{W}_k$ and $\mathcal{V}_k$, respectively.
The shape matrix of the ellipsoidal sets $\mathcal{W}_k$ and $\mathcal{V}_k$ are
\begin{eqnarray}
\nonumber\mQ_k&=&10\left[
                     \begin{array}{cccc}
                       \frac{T^3}{3} & 0 & \frac{T^2}{2} & 0 \\
                       0 & \frac{T^3}{3} & 0 &  \frac{T^2}{2} \\
                      \frac{T^2}{2} & 0 & T & 0 \\
                       0 & \frac{T^2}{2} & 0 & T \\
                     \end{array}
                   \right]
\\
\nonumber\mR_k&=&\left[
                   \begin{array}{cc}
                     100 & 0 \\
                     0 & 0.5 \\
                   \end{array}
                 \right].
\end{eqnarray}
In the simulation, the noise obeys the unform distribution, see Fig. \ref{fig_04}. Assume the initial estimation shape matrix is $\mP_0=200\mI$, and the initial estimation state $\hat{\vx}_0$ is a random disturbance around the true state $\vx_0$.

The simulation results in this example include three parts: the first part is about the robustness, the second part is about the size of the estimated  ellipsoid, and the third part is about the running time.

\begin{figure}[ht]
\vbox to 5cm{\vfill \hbox to \hsize{\hfill
\scalebox{0.5}[0.5]{\includegraphics{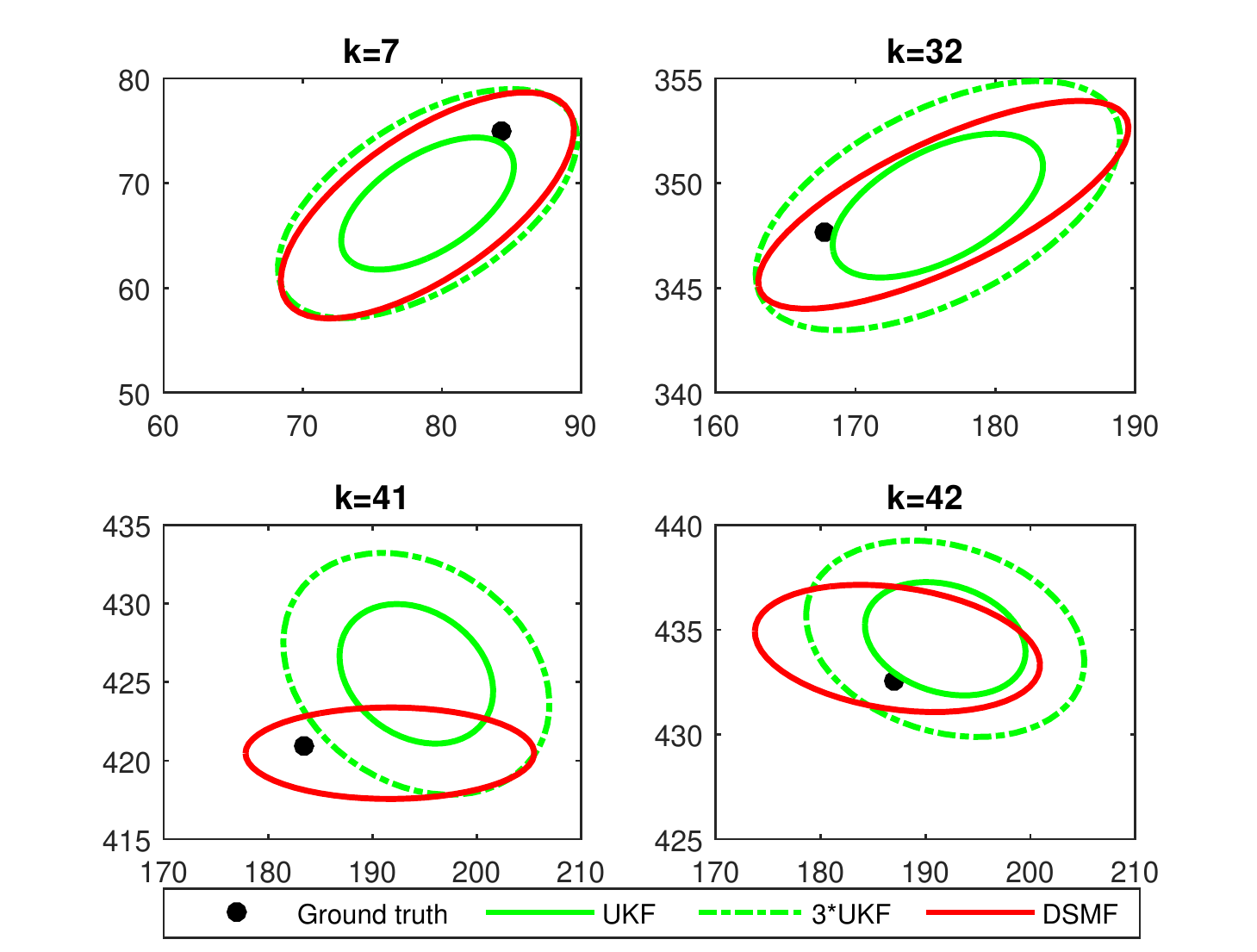}} \hfill}\vfill}
\caption{The estimated ellipsoid by DSMF and UKF at different time step.}\label{fig_05}
\end{figure}

\begin{figure}[ht]
\vbox to 5cm{\vfill \hbox to \hsize{\hfill
\scalebox{0.5}[0.5]{\includegraphics{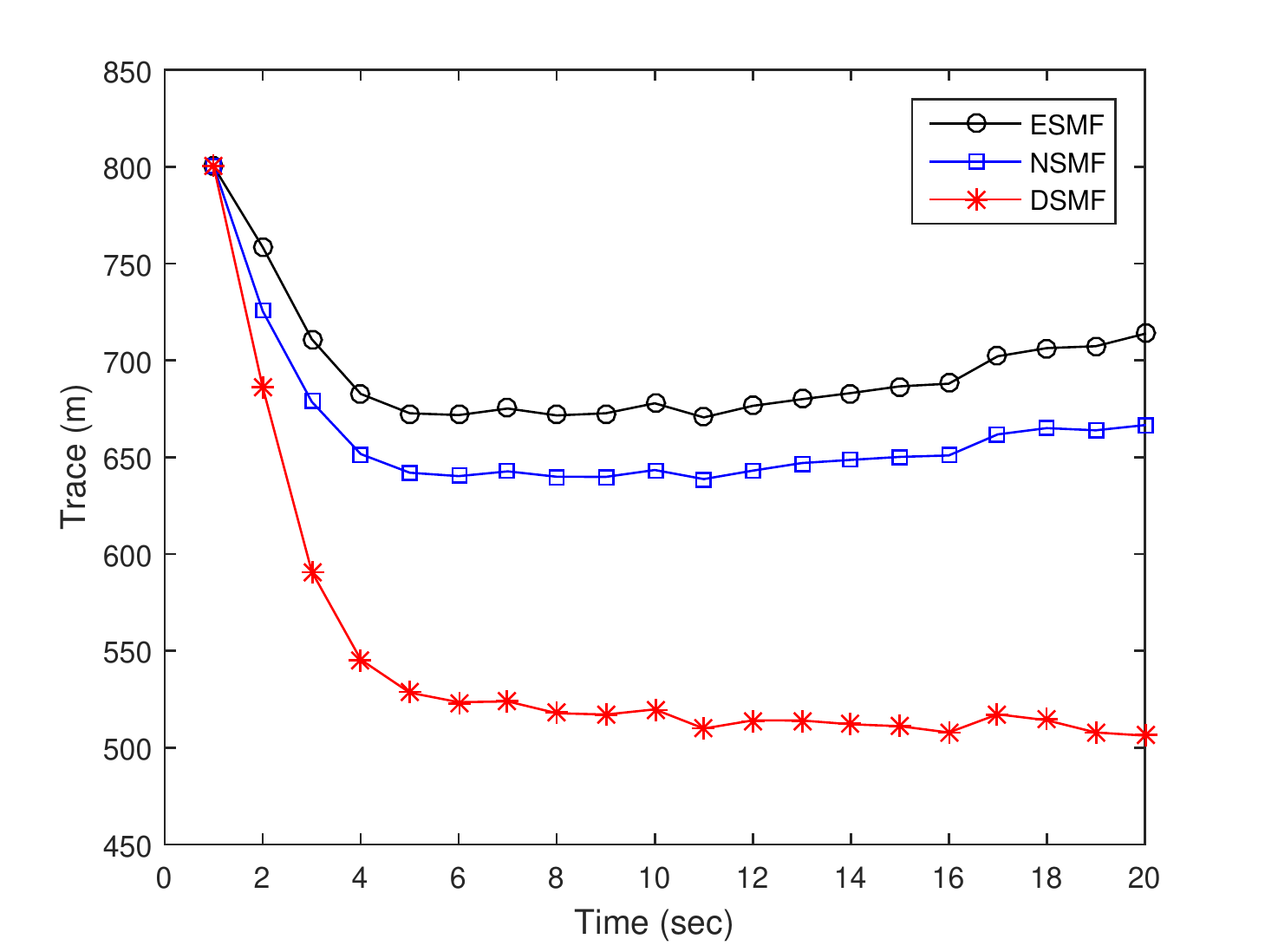}} \hfill}\vfill}
\caption{The size of the estimated ellipsoid by the different methods.}\label{fig_06}
\end{figure}

\begin{figure}[ht]
\vbox to 5cm{\vfill \hbox to \hsize{\hfill
\scalebox{0.5}[0.5]{\includegraphics{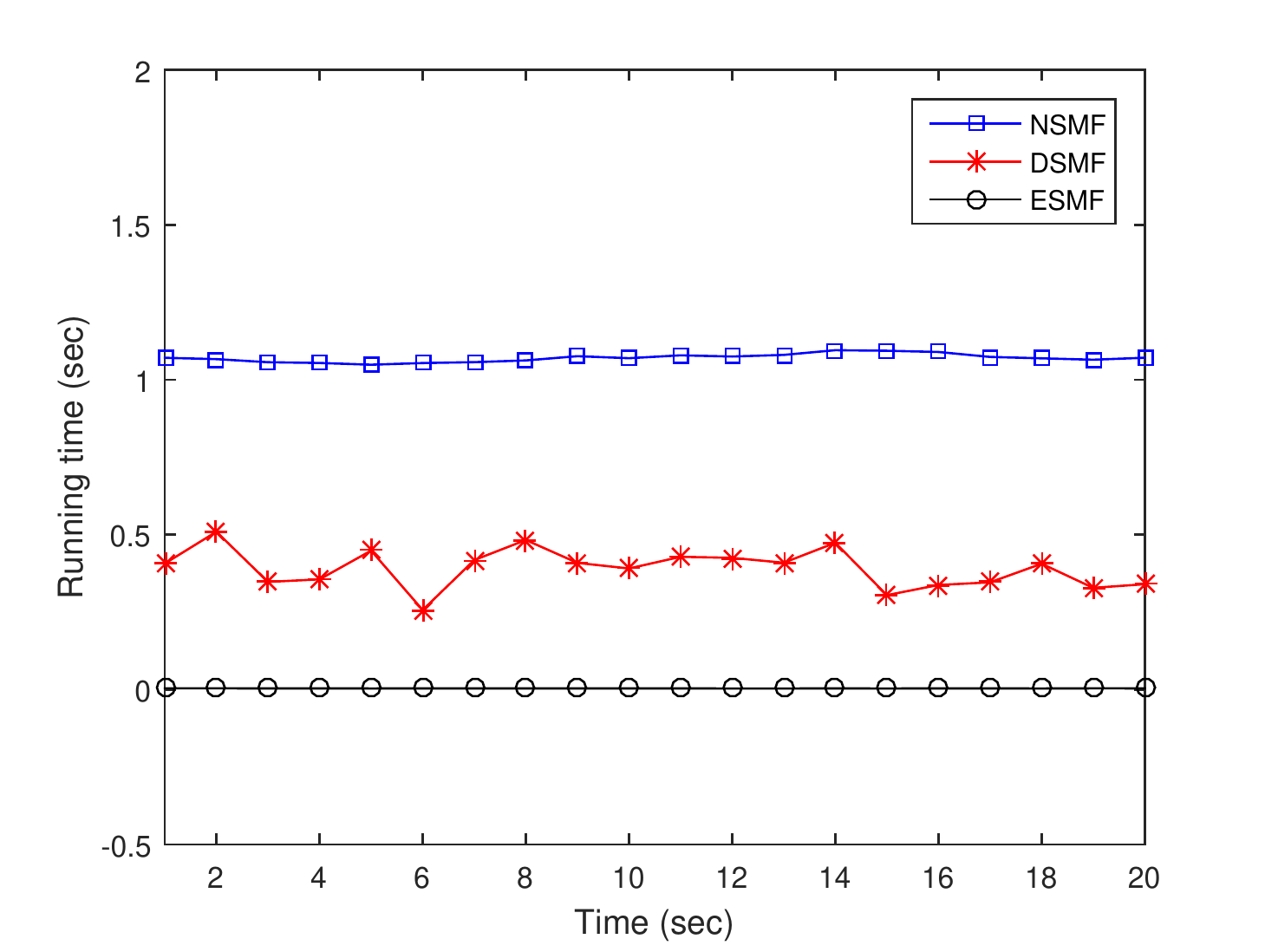}} \hfill}\vfill}
\caption{The running time by the different methods.}\label{fig_07}
\end{figure}

\begin{itemize}
  \item Fig. \ref{fig_05} shows the estimated ellipsoid by DSMF and UKF at time steps 7, 32, 41, 42. the dotted ellipsoid indicates the three times the confidence ellipsoid by UKF. The (stochastic) confidence ellipsoid provided by the UKF is indeed tighter than their deterministic counterparts computed via DSMF, but it does not guarantee the containment of the true state, even if the three times of the confidence ellipsoid by UKF.
  \item In Fig. \ref{fig_06}, the size $tr(\mP_{k})$ of the estimated ellipsoid $\mathcal{E}_{k}$ is plotted as a function of the time step. It shows that the size of the estimated ellipsoid can quickly converge to a stable value, and DSMF has the smallest estimated ellipsoid. The reason may be that DSMF fully extracts the  properties of the nonlinear function by solving a semi-infinite optimization problem, rather than linearizing it.
  \item The running times of the different set membership filters are plotted in Fig. \ref{fig_07}, respectively. From Fig. \ref{fig_07}, it shows that the running time of DSMF stays between  ESMF and NSMF. The reason is that ESMF does not need to solve the optimization problem, and DSMF and NSMF need to solve some optimization problems by the first order method and SDP, respectively. In order to get a tradeoff between the performance and running time, DSMF method is the best choice.
\end{itemize}


\subsection{Mobile Robot Localization}
Consider the localization of a mobile robot operating
in planar environments. Then, the nonlinear state function $f(\vx_k)$ defined in (\ref{equ_5_9}) is given as
\begin{align}
\nonumber f_k(\vx_k)=\left[
                       \begin{array}{c}
                         p_k^x-\frac{u_p}{u_r}(\sin(\theta_k)-\sin(\theta_k+T_0u_r)) \\
                         p_k^y+\frac{u_p}{u_r}(\cos(\theta_k)-\cos(\theta_k+T_0u_r)) \\
                         \theta_k+T_0u_r \\
                       \end{array}
                     \right],
\end{align}
where $u_p=0.085$ and $u_r=0.015$ are the motion command to control the translational velocity and  the rotational velocity, respectively. Here, $T_0$ is the sampling period with $T_0=1$.  The state noise $\vw_k$ in (\ref{equ_5_9})  belongs to an ellipsoid $\mathcal{W}_k$, and its shape matrix $\mQ_k=\diag(10e-7,10e-7, 10e-8)$. The measurement equation is defined in (\ref{equ_5_10}) and the landmark location is $s^x=50,s^y=50$. Assume the measurement noise $\vv_k$ is bounded by an ellipsoid $\mathcal{V}_k$, and its shape matrix $\mR_k=\diag(1,1)$. In our simulation, the initial state $\vx_0$ is $[10~ 10~ 1]^T$. Moreover, the initial target state estimate $\hat{\vx}_0$ is a combination of the actual target state and a random bias vector, and initial shape matrix $\mP_0=\diag(1,1,0.1)$. The root mean square error (RMSE) is defined as follows
\begin{align}
\textmd{RMSE}_k=\sqrt{\frac{\sum_{i=1}^L(\hat{x}_k^i-x_k^i)^2}{L}},
\end{align}
where $x_k^i$ is the true state and $\hat{x}_k^i$ is the state estimation at $i$th Monte Carlo, in addition, $L=50$ in this simulation.

\begin{figure}[ht]
\vbox to 5cm{\vfill \hbox to \hsize{\hfill
\scalebox{0.5}[0.5]{\includegraphics{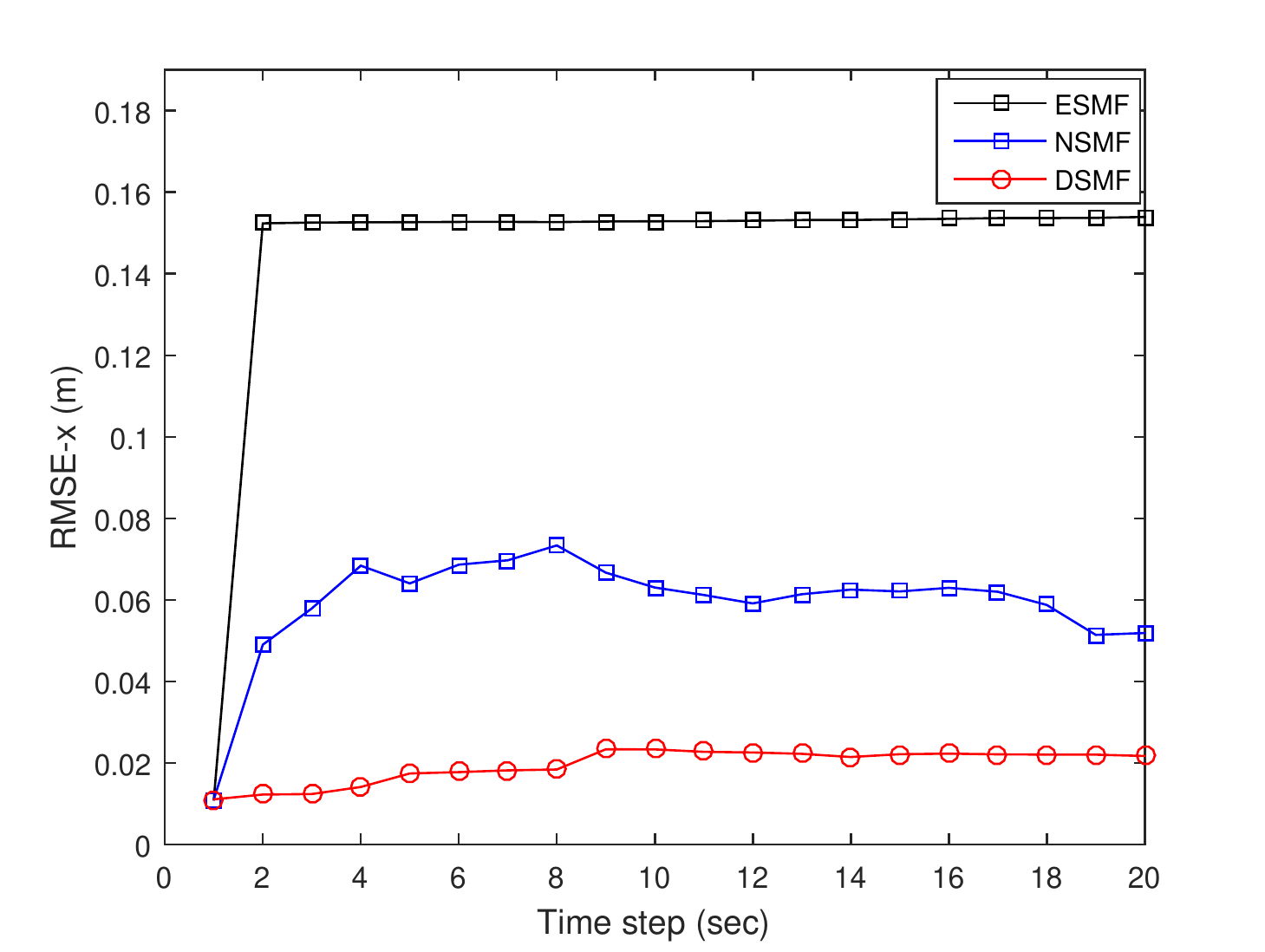}} \hfill}\vfill}
\caption{The RMSE of the state $p_k^x$ estimation is plotted as a function of time
step.}\label{fig_10}
\end{figure}

\begin{figure}[ht]
\vbox to 5cm{\vfill \hbox to \hsize{\hfill
\scalebox{0.5}[0.5]{\includegraphics{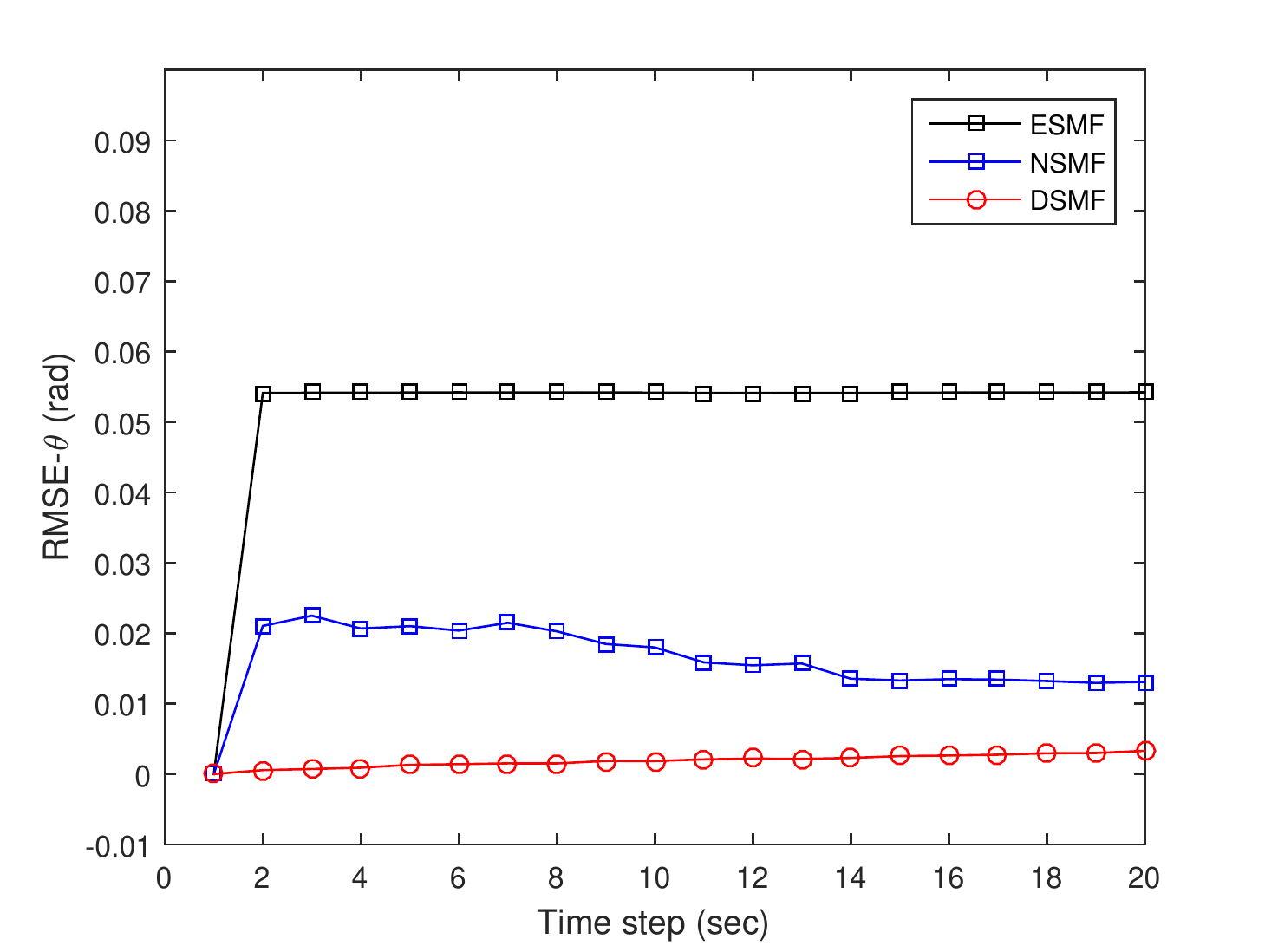}} \hfill}\vfill}
\caption{The RMSE of the state $\theta_k$ estimation is plotted as a function of time
step. }\label{fig_11}
\end{figure}

In this example, the RMSE of the state estimation along $x$ direction and $\theta$ direction is plotted as a function of the time steps in Figs. \ref{fig_10}-\ref{fig_11}. They show that the RMSE of DSMF is less than that of ESMF, NSMF. The reason may be that DSMF uses the information of nonlinear functions sufficiently, and derives a tighter ellipsoid to cover the state.
\section{Conclusion}\label{sec_7}
In this paper, we have proposed the dual set membership filter for nonlinear dynamic system with unknown but bounded noise. It involves two steps, prediction step and measurement updated step, to recursively compute a bounding ellipsoid to cover the true state. We  translated the nonlinear system  into the linear system  by solving the semi-infinite programming, rather than linearizing the nonlinear function. Moreover, the semi-infinite programming can be solved efficiently by the Frank-Wolfe method, then it is suitable for large scale problems, and the low computing complexity in each step  makes it working on-line. In addition,  the proposed filter has bee applied to mobile robot localization. Finally, three illustrative examples in the simulations show that the proposed filter performs better. Possible research direction may include multi-target tracking, sensor management  in the sensor networks based on the proposed filter.
\section{Appendix}
\textbf{The proof of Lemma \ref{lem_1}}:
\begin{proof}
Defining an operator $G:\mathbb{R}^{d\times d}\to C(\Omega)$ by
\begin{eqnarray}
\nonumber [G(\mH)](\tilde{\vy})\triangleq \tilde{\vy}^T\mH\tilde{\vy}-d ~on ~\Omega,
\end{eqnarray}
where $C(\Omega)$ denotes the Banach space of all continuous on $\Omega$. Thus we can write problem (\ref{equ_63})-(\ref{equ_64}) as
\begin{eqnarray}
\label{equ_65}&&\min ~-\mbox{logdet}(\mH)\\
\label{equ_66}&& \mbox{s.t.}~ G(\mH)\leq 0.
\end{eqnarray}
Let us introduce the Lagrangian $L:\mathbb{R}^{d\times d}\times M(\Omega)\to \mathbb{R}$ for problem (\ref{equ_65})-(\ref{equ_66}) as
\begin{eqnarray}
\label{equ_67} L(\mH,\mu)=-\mbox{logdet}(\mH)+\langle \mu, G(\mH) \rangle,
\end{eqnarray}
where $M(\Omega)$ is the dual space of $C(\Omega)$, and  $\langle \cdot, \cdot \rangle$ denotes a duality pair between $C(\Omega)$ and $M(\Omega)$  as
\begin{eqnarray}
\nonumber \langle \mu, \nu \rangle =\int_{\Omega}\nu(\tilde{\vy})d\mu.
\end{eqnarray}
Let us denote the dual operator of the derivative $\nabla G(\mH)$ \cite{Ito-Liu-Teo00}, which is given by
\begin{eqnarray}
\label{equ_68} \langle \nabla \mu,G(\mH) \rangle=\int_{\Omega}\tilde{\vy}\tilde{\vy}^T d\mu.
\end{eqnarray}
By KKT conditions, we know that
\begin{eqnarray}
\label{equ_74}\nabla_{H} L(\mH^{*},\mu^{*})=-{H^{*}}^{-1}+\langle\nabla \mu^{*},G(\mH^{*})\rangle=0\\
\label{equ_75}\langle \mu^{*}, \tilde{\vy}^T\mH^{*}\tilde{\vy}-d \rangle=0, ~\mu^{*}\geq 0.
\end{eqnarray}
Substituting (\ref{equ_68}) into (\ref{equ_74}), we obtain the final result.
\end{proof}
\textbf{The proof of Lemma \ref{lem_2}}:
\begin{proof}
According to definition of Lagrangian function and its gradient in (\ref{equ_67}) and  (\ref{equ_74}), we conclude that
\begin{eqnarray}
\label{equ_76}\min_{\mH}L(\mH,\mu)=\mbox{logdet}\Big(\int_{\Omega}\tilde{\vy}\tilde{\vy}^T d\mu\Big)+\int_{\Omega}\Big(\tilde{\vy}^T(\int_{\Omega}\tilde{\vy}\tilde{\vy}^Td\mu)^{-1}\tilde{\vy}-d\Big) d\mu.
\end{eqnarray}
Since
\begin{eqnarray}
\nonumber\int_{\Omega}\tilde{\vy}^T\Big(\int_{\Omega}\tilde{\vy}\tilde{\vy}^Td\mu\Big)^{-1}\tilde{\vy} d\mu&=&\int_{\Omega}tr\Big(\tilde{\vy}^T(\int_{\Omega}\tilde{\vy}\tilde{\vy}^Td\mu)^{-1}\tilde{\vy}\Big) d\mu\\
\nonumber&=&\int_{\Omega}tr\Big((\int_{\Omega}\tilde{\vy}\tilde{\vy}^Td\mu)^{-1}\tilde{\vy}\tilde{\vy}^T\Big)d\mu\\
\nonumber&=&tr\Big((\int_{\Omega}\tilde{\vy}\tilde{\vy}^Td\mu)^{-1}\int_{\Omega}\tilde{\vy}\tilde{\vy}^Td\mu\Big)\\
\nonumber&=&d,
\end{eqnarray}
we can rewrite (\ref{equ_76}) as
\begin{eqnarray}
\min_{\mH}L(\mH,\mu)=\mbox{logdet}\Big(\int_{\Omega}\tilde{\vy}\tilde{\vy}^T d\mu\Big)+d-d\int_{\Omega}d\mu.
\end{eqnarray}
This leads to a first dual problem
\begin{eqnarray}
&&\max_{\mu}~\mbox{logdet}\Big(\int_{\Omega}\tilde{\vy}\tilde{\vy}^T d\mu\Big)+d-d\int_{\Omega}d\mu\\
&&\mbox{s.t.}~\mu\geq 0.
\end{eqnarray}
In fact, we can restrict $\mu$ to satisfy $\int_{\Omega}d\mu=1$, then we obtain the final result.
\end{proof}
\textbf{The proof of Proposition \ref{pro_1}}:
\begin{proof}
Let $\mH_{y}=\mP^{-1}$, then the optimization problem (\ref{equ_61})-(\ref{equ_62}) can be rewritten as
\begin{eqnarray}
\label{equ_77} &&\min ~-\mbox{logdet}(\mH_{y})\\
\label{equ_78}  && \mbox{s.t.}~ (\vy-\hat{\vy})^T\mH_y(\vy-\hat{\vy})\leq n,
 ~\forall~ \vy\in\mathcal{C}.
\end{eqnarray}
The constrain (\ref{equ_78}) is equivalent to
\begin{align}
\nonumber\left[
  \begin{array}{c}
    \vy \\
    1 \\
  \end{array}
\right]^T\left[
         \begin{array}{cc}
           \mI & 0 \\
           -\hat{\vy}^T & 1 \\
         \end{array}
       \right]\left[
                \begin{array}{cc}
                  \mH_y & 0 \\
                  0 & 1 \\
                \end{array}
              \right]\left[
         \begin{array}{cc}
           \mI & -\hat{\vy} \\
            0& 1 \\
         \end{array}
       \right]\left[
  \begin{array}{c}
    \vy \\
    1 \\
  \end{array}
\right]\leq d,~\forall~ \vy\in\mathcal{C}
\end{align}
where $d=n+1$. If we define $\tilde{y}=[y^T ~1]^T$, it shows that there exists an ellipsoid with shape matrix $\tilde{\mH}$,
\begin{align}
\tilde{\mH}:=\left[
         \begin{array}{cc}
           \mI & 0 \\
           -\hat{\vy}^T & 1 \\
         \end{array}
       \right]\left[
                \begin{array}{cc}
                  \mH_y & 0 \\
                  0 & 1 \\
                \end{array}
              \right]\left[
         \begin{array}{cc}
           \mI & -\hat{\vy} \\
            0& 1 \\
         \end{array}
       \right],
\end{align}
which covers the set $\mathcal{C}\times\{1\}$. Note that $-\mbox{logdet}(\mH_y)=-\mbox{logdet}(\tilde{\mH})\geq-\mbox{logdet}(\tilde{\mH}^{*})$, where $\tilde{\mH}^{*}$ is the optimal solution of problem
\begin{eqnarray}
\label{equ_79} &&\min ~-\mbox{logdet}(\tilde{\mH})\\
\label{equ_80}  && \mbox{s.t.}~ \tilde{\vy}^T\tilde{\mH}\tilde{\vy}\leq d,
 ~\forall~ \tilde{\vy}\in\mathcal{C}\times\{1\}.
\end{eqnarray}
From Lemma \ref{lem_1} and Lemma \ref{lem_2}, we know
\begin{align}
\nonumber \tilde{\mH}^{*}&=\Big(\int \tilde{\vy} \tilde{\vy}^Td\mu^{*}\Big)^{-1}=\left[
                                                    \begin{array}{cc}
                                                      \int \vy \vy^Td\mu^{*}  &  \int \vy d\mu^{*} \\
                                                       \int \vy^T  d\mu^{*} &  1 \\
                                                    \end{array}
                                                  \right]^{-1}\\
 \label{equ_81}&=\left[
         \begin{array}{cc}
           \mI & 0 \\
           -\int \vy^T  d\mu^{*} & 1 \\
         \end{array}
       \right]\left[
                \begin{array}{cc}
                 (\int \vy \vy^Td\mu^{*}-\hat{\vy}^{*}(\hat{\vy}^{*})^T)^{-1} & 0 \\
                  0 & 1 \\
                \end{array}
              \right]\left[
         \begin{array}{cc}
           \mI & -\int \vy  d\mu^{*} \\
            0& 1 \\
         \end{array}
       \right]
\end{align}
Let us set $\mH_y^{*}=(\int \vy \vy^Td\mu^{*}-\hat{\vy}^{*}(\hat{\vy}^{*})^T)^{-1}$ and $\hat{\vy}^{*}=\int \vy  d\mu^{*}$, we note that this leads that
$-\mbox{logdet}(\mH_y^{*})=-\mbox{logdet}(\tilde{\mH}^{*})$, so that
\begin{align}
 \label{equ_83}-\mbox{logdet}(\mH_y)\geq -\mbox{logdet}(\mH_y^{*}).
\end{align}
Since $\mH_y^{*}$ satisfies the constraint (\ref{equ_80}), it leads to
\begin{align}
 \label{equ_82}\tilde{\vy}^T\tilde{\mH}^{*}\tilde{\vy}\leq d.
\end{align}
Substituting (\ref{equ_81}) into (\ref{equ_82}) yields
\begin{align}
(\vy-\hat{\vy}^{*})^T\mH_y^{*}(\vy-\hat{\vy}^{*})\leq n,
 ~\forall~ \vy\in\mathcal{C},
\end{align}
so that the ellipsoid with center $\hat{\vy}^{*}$ and shape matrix $\mH_y^{*}$ contains the set $\mathcal{C}$. According to (\ref{equ_83}), it proves the minimality of this ellipsoid.
\end{proof}

\end{document}